\documentclass{amsart}

\usepackage{amsfonts,amsmath,amssymb,mathrsfs,verbatim}

\numberwithin{equation}{section}

\newtheorem{lemma}{Lemma}[section]

\newtheorem{theorem}{Theorem}[section]

\newtheorem{proposition}{Proposition}[section]
\newtheorem{conjecture}{Conjecture}[section]

\newcommand{\Z}{\mathbb Z}
\newcommand{\abs}[1]{\lvert#1\rvert}
\newcommand{\set}[1]{[#1]}
\newcommand{\Zet}[1]{\Z_{#1}}
\newcommand{\pqinf}{(p;p)_{\infty}(q;q)_{\infty}}
\newcommand{\pq}[1]{(p;p)_{\infty}^{#1}(q;q)_{\infty}^{#1}}
\renewcommand{\d}[1]{\frac{\textup{d}#1}{#1}}
\renewcommand{\Re}{\textup{Re}}
\newcommand{\e}{\textup{e}}
\newcommand{\Complex}{\mathbb C}

\newcommand{\R}{\mathbb R}
\newcommand{\la}{\lambda}
\newcommand{\im}{\hspace{0.5pt}\textup{i}\hspace{0.5pt}}
\newcommand{\Res}{\textup{Res}}
\newcommand{\An}{\mathbb A}
\newcommand{\F}{\mathcal F}
\newcommand{\A}{\mathscr A}
\newcommand{\C}{\mathscr C}
\newcommand{\varw}{\textup w}
\newcommand{\T}{\mathbb T}

\begin{document}

\title[Elliptic beta integrals]
{Inversions of integral operators and elliptic beta integrals
on root systems}

\author{Vyacheslav P. Spiridonov}\thanks{V.P.S. is
supported in part by the Russian Foundation for Basic Research
(RFBR) grant no. 03-01-00781; S.O.W. is supported by the Australian
Research Council}
\address{Bogoliubov Laboratory of Theoretical Physics,
JINR, Dubna, Moscow Region 141980, Russia and
Max--Planck--Institut f\"ur
Mathematik, Vivatsgasse 7, D-53111, Bonn, Germany}
\email{spiridon@thsun1.jinr.ru}
\author{S. Ole Warnaar}
\address{Department of Mathematics and Statistics,
The University of Melbourne, VIC 3010, Australia}
\email{warnaar@ms.unimelb.edu.au}

\date{September 2004}



\begin{abstract}
We prove a novel type of inversion formula for elliptic hypergeometric
integrals associated to a pair of root systems. Using the (A,C)
inversion formula to invert one of the known C-type elliptic beta integrals,
we obtain a new elliptic beta integral for the root system of type A. Validity
of this integral is established by a different method as well.
\end{abstract}

\maketitle

\tableofcontents

\section{Introduction}

Beta-type integrals are fundamental objects of applied analysis, with
numerous applications in pure mathematics and mathematical physics.
The classical Euler beta integral
\begin{equation*}
\int_0^1 t^{x-1}(1-t)^{y-1}\textup{d}t=\frac{\Gamma(x)\Gamma(y)}{\Gamma(x+y)},
\qquad \min\{\Re(x),\Re(y)\}>0,
\end{equation*}
determines the measure for the Jacobi family of orthogonal
polynomials expressed as certain $_2F_1$ hypergeometric functions \cite{AAR}.
Its multi-dimensional extension due to Selberg \cite{Selberg44}
plays an important role in harmonic analysis on root systems,
the theory of special functions of many variables,
the theory of random matrices, and so forth.

Important generalizations of beta integrals arise in the
theory of basic or $q$-hypergeometric functions.
The Askey--Wilson $q$-beta integral depends on four independent parameters and a
base $q$, and fixes the orthogonality measure for the Askey--Wilson
polynomials, the most general family of classical single-variable orthogonal
polynomials \cite{AW85}.
Closely related to the Askey--Wilson integral is the
integral representation for a very-well-poised $_8\phi_7$ basic hypergeometric
series found by Nassrallah and Rahman \cite{NR85}. Through specialization this
led Rahman to the discovery of a one-parameter extension of the Askey--Wilson
integral \cite{Rahman86}.
Finally, several multi-dimensional generalizations of the
Askey--Wilson and Rahman integrals, including a $q$-Selberg integral,
were found by Gustafson \cite{Gustafson90,Gustafson92,Gustafson94}.
For some time, these multi-dimensional $q$-beta integrals were believed to be
the most general integrals of beta type.

A new development in the field was initiated by the first author with the discovery
of an elliptic generalization of Rahman's $q$-beta integral \cite{Spiridonov01}.
This elliptic beta integral depends on five free parameters and two basic variables
--- or elliptic moduli --- $p$ and $q$.
As a further development two $n$-dimensional elliptic beta integrals
associated to the C$_n$ root system were proposed by van Diejen and the
first author \cite{vDS00,vDS01}. In the $p\to 0$ limit these integrals reduce to
Gustafson's C$_n$ $q$-beta integrals.
More elliptic beta integrals, all related to either the A$_n$ or C$_n$ root systems
and all but one generalizing integrals of Gustafson \cite{Gustafson92,Gustafson94}
and Gustafson and Rakha \cite{GR00},
were subsequently given in \cite{Spiridonov04}.

Roughly, $n$-dimensional elliptic beta integrals come in three different types.
Most fundamental are the type-I integrals. These contain $2n+3$ free parameters
(as well as the bases $p$ and $q$), and one of the A$_n$ integrals of
\cite{Spiridonov03b} and one of the C$_n$ integrals of \cite{vDS01} are of type I.
The first complete proofs were found by Rains \cite{Rains03} who derived them
as a consequence of a symmetry transformation for more general elliptic
hypergeometric integrals.
More elementary proofs using difference equations
were subsequently given in \cite{Spiridonov04}.
The elliptic beta integrals of type II contain less than $2n+3$
parameters and can be deduced from type I integrals via the composition
of higher-dimensional integrals
\cite{vDS01,Gustafson90,Gustafson94,GR00,Spiridonov03b}.
The second C$_n$ elliptic beta integral of
\cite{vDS01} (see also \cite{vDS00}), depending on six parameters
(only $5$ when $n=1$), provides an example of a type II integral.
Finally, type III elliptic beta integrals arise through the
computation of $n$-dimensional determinants with entries composed of
one-dimensional integrals \cite{Spiridonov03b}.
Originally, all of the above beta integrals were defined
for bases $p$ and $q$ inside the unit circle (due to the use of
the standard elliptic gamma function described in the next section).
Another class of elliptic hypergeometric integrals, which
are well defined in the larger region $\abs{p}<1$, $\abs{q}\leq 1$
(by employing a different elliptic gamma function),
has been introduced in \cite{Spiridonov03b}.
We shall not discuss here the corresponding elliptic beta
integrals, and refer the reader to \cite{vDS03,Spiridonov04} for more details.

Further progress on the subject is associated with symmetry transformations
of elliptic  hypergeometric integrals. Certain hypergeometric identities are
well-known to be related to the notion of matrix inversions and Bailey pairs.
At the level of hypergeometric series --- ordinary, basic or elliptic ---
the Bailey pair machinery allows for the derivation of infinite sequences
of symmetry transformations \cite{Andrews84,Andrews01,Spiridonov02,W01,W03}.
A formulation of the notion of Bailey pairs for integrals was proposed in
\cite{Spiridonov04a} (on the basis of a transformation for univariate elliptic
hypergeometric integrals proved in \cite{Spiridonov03b}). Using the univariate
elliptic beta integral, this led to a binary tree of identities for multiple
elliptic hypergeometric integrals. A generalization of these results to elliptic
hypergeometric integrals labelled by root systems has been one of the motivations
for the present paper. Indeed, the integral analogues of the matrix inversions
underlying the Bailey transform for series are provided by the integral inversions
of this paper.

A powerful set of symmetry transformations relating elliptic hypergeometric
integrals of various dimensions was introduced by Rains \cite{Rains03}.
He proved the latter in an elegant manner by reducing the problem to
determinant evaluations on a dense set of parameters. Although some of the
Rains transformations can be reproduced with the help of the Bailey type
technique, a complete correspondence between these two sets of
identitites has not been established yet.

More specifically, we provide the following new framework for viewing
elliptic beta integrals on root systems.
First, we introduce certain multi-dimensional integral transformations
with integration kernels determined by the structure of the type I elliptic beta
integrals on the A$_n$ and C$_n$ root systems.
The A$_n$ and C$_n$ elliptic beta integrals then acquire the new interpretation
as examples for which these integral transformations can be performed explicitly.
Second, we prove two theorems describing inversions of the corresponding
integral operators on a certain class of functions and conjecture a third
inversion formula. These inversion formulas naturally carry two root system labels,
our three results corresponding to the pairs (A$_n$,A$_n$), (A$_n$,C$_n$)
and (C$_n$,A$_n$). Third, using the (A$_n$,C$_n$) inversion formula we `invert'
the type I C$_n$ beta integral to prove a new type I A$_n$ elliptic beta integral.
It appears that this exact integration formula is new even at the
$q$-hypergeometric and plain hypergeometric levels.
Finally, for completeness, we give
an alternative proof of this integral using the method for proving
type I integrals developed in \cite{Spiridonov04}.

In the univariate case, all three integral inversions coincide and the
resulting formula establishes the inversion of the integral Bailey
transform of \cite{Spiridonov04a}.
Also our multi-variable integral transformations
on root systems can be put into the framework of integral Bailey pairs.
This will be the topic of a subsequent publication together with a
consideration of integral operators associated with the type II
elliptic beta integrals.

\section{Notation and preliminaries}

Throughout this paper $p,q\in\Complex$ such that
\begin{equation}\label{Mdef}
M:=\max\{\abs{p},\abs{q}\}<1.
\end{equation}

For fixed $p$ and $q$, and $z\in\Complex\backslash \{0\}$ the elliptic gamma
function is defined as \cite{Ruijsenaars97}
\begin{equation}\label{ellGamma}
\Gamma(z;p,q)=\prod_{\mu,\nu=0}^{\infty}
\frac{1-z^{-1}p^{\mu+1}q^{\nu+1}}{1-zp^{\mu} q^{\nu}},
\end{equation}
and satisfies
\begin{subequations}
\begin{align}
\Gamma(z;p,q)&=\Gamma(z;q,p), \\
\Gamma(z;p,q)&=\frac{1}{\Gamma(pq/z;p,q)}.
\label{Gammasymm}
\end{align}
\end{subequations}
Defining the theta function
\begin{equation*}
\theta(z;p)=(z,p/z;p)_{\infty},
\end{equation*}
where
\begin{equation*}
(a;q)_{\infty}=\prod_{k=0}^{\infty}(1-aq^k)
\quad\text{and}\quad
(a_1,\dots,a_k;q)_n=(a_1;q)_n\cdots (a_k;q)_n,
\end{equation*}
it follows that
\begin{equation}\label{Gammaqd}
\Gamma(qz;p,q)=\theta(z;p)\Gamma(z;p,q), \quad
\Gamma(pz;p,q)=\theta(z;q)\Gamma(z;p,q)
\end{equation}
and
\begin{equation}\label{reflex}
\Gamma(z;p,q)\Gamma(z^{-1};p,q)=\frac{1}{\theta(z;p)\theta(z^{-1};q)}.
\end{equation}
A useful formula needed repeatedly for calculating residues is
\begin{equation}\label{Glim}
\lim_{z\to a}(1-z/a)\Gamma(z/a;p,q)=
\frac{1}{\pqinf}.
\end{equation}

For $n$ an integer the elliptic shifted factorial is defined by \cite{W02}
\begin{equation*}
(a;q,p)_n=\frac{\Gamma(aq^n;p,q)}{\Gamma(a;p,q)}
\end{equation*}
(this was denoted as $\theta(a;p;q)_n$ in \cite{vDS00,vDS01,Spiridonov03b}).
When $n$ is non-negative it may also by written as
\begin{equation*}
(a;q,p)_n=\prod_{j=0}^{n-1}\theta(aq^j;p).
\end{equation*}
For both the elliptic gamma function and the elliptic shifted
factorial we employ standard condensed notation, i.e.,
\begin{align*}
\Gamma(z_1,\dots,z_k;p,q)&=\Gamma(z_1;p,q)\cdots\Gamma(z_k;p,q) \\
(a_1,\dots,a_k;q,p)_n&=(a_1;q,p)_n\cdots(a_k;q,p)_n.
\end{align*}
Also, we will often suppress the $p$ and $q$ dependence
and write $\Gamma(z)=\Gamma(z;p,q)$,
$\theta(z)=\theta(z;p)$ and $(a)_n=(a;q,p)_n$.

As a final notational point we write the sets $\{1,\dots,n\}$ and
$\{0,\dots,n-1\}$ as $\set{n}$ and $\Zet{n}$, and
adopt the convention that $\mu$ and $\nu$ are
non-negative integers.

\section{Elliptic beta integrals}\label{secbetaint}

The single-variable elliptic beta integral
--- due to the first author~\cite{Spiridonov01} ---
corresponds to the following generalization of the celebrated
Rahman integral \cite{Rahman86} (obtained as an important
special case of the Nassrallah--Rahman integral \cite{NR85}).
Let $t_1,\dots,t_6\in\Complex$ such that $t_1\cdots t_6=pq$ and
\begin{equation}\label{cond}
\max\{\abs{t_1},\dots,\abs{t_6}\}<1,
\end{equation}
and let $\T$ denote the positively oriented unit circle.
Then
\begin{equation}\label{elbetasymm}
\frac{1}{2\pi\im}
\int_{\T}\frac{\prod_{i=1}^6 \Gamma(t_iz,t_iz^{-1})}
{\Gamma(z^2,z^{-2})}\d{z}
=\frac{2\prod_{1\leq i<j\leq 6}\Gamma(t_it_j)}{\pqinf}.
\end{equation}
Defining $T=t_1\cdots t_5$, eliminating
$t_6$ and using the symmetry \eqref{Gammasymm},
this may also be put in the form
\begin{equation}\label{elbeta}
\frac{1}{2\pi\im}
\int_{\T}\frac{\prod_{i=1}^5 \Gamma(t_iz,t_iz^{-1})}
{\Gamma(z^2,z^{-2},Tz,Tz^{-1})}\d{z}
=\frac{2\prod_{1\leq i<j\leq 5}\Gamma(t_it_j)}
{\pqinf\prod_{i=1}^5\Gamma(Tt_i^{-1})}
\end{equation}
from which the Rahman integral follows by
letting $p$ (or, equivalently, $q$) tend to $0$.
For \eqref{elbeta} to be valid we must of course
replace \eqref{cond} by
\begin{equation*}
\max\{\abs{t_1},\dots,\abs{t_5},\abs{pq/T}\}<1.
\end{equation*}

Two multivariable generalizations of \eqref{elbetasymm}
associated with the root systems of type A and C will be needed.
In order to state these we require some further notation.
Throughout, $n$ will be a fixed positive integer,
$z=(z_1,\dots,z_n)$ and
\begin{equation*}
\d{z}=\d{z_1}\cdots\d{z_n}.
\end{equation*}
Whenever the variable $z_{n+1}$ occurs it will be fixed by
$z_1\cdots z_{n+1}=1$ unless stated otherwise.
For reasons of printing economy we also employ the notation $f(z_i^{\pm})$ for
$f(z_i,z_i^{-1})$, $f(z_i^{\pm}z_j^{\pm})$ for
$f(z_1z_j,z_iz_j^{-1},z_i^{-1}z_j,z_i^{-1}z_j^{-1})$ and so on.

The A$_n$ generalization of \eqref{elbetasymm} depends on
$2n+4$ complex parameters $t_1,\dots,t_{n+2}$ and $s_1,\dots,s_{n+2}$
such that $\max\{\abs{t_1},\dots,\abs{t_{n+2}},\abs{s_1},\dots,\abs{s_{n+2}}\}<1$
and $ST=pq$ for $T=t_1\cdots t_{n+2}$ and $S=s_1\cdots s_{n+2}$.
Hence we effectively have only $2n+3$ free parameters, making it an
elliptic beta integral of type I;
\begin{multline}\label{Anbeta}
\frac{1}{(2\pi\im)^n}\int_{\T^n}
\prod_{i=1}^{n+2}\prod_{j=1}^{n+1}\Gamma(s_iz_j,t_iz_j^{-1})
\prod_{1\leq i<j\leq n+1}\frac{1}{\Gamma(z_iz_j^{-1},z_i^{-1}z_j)}
\d{z} \\
=\frac{(n+1)!}{\pq{n}}
\prod_{i=1}^{n+2}\Gamma(Ss_i^{-1},Tt_i^{-1})
\prod_{i,j=1}^{n+2}\Gamma(s_it_j).
\end{multline}
As already mentioned in the introduction this integral was conjectured by
the first author \cite{Spiridonov03b}
and subsequently proven by Rains \cite[Corollary 4.2]{Rains03}
and by the first author \cite[Theorem 3]{Spiridonov04}.

The type I elliptic beta integral for the root system C$_n$
depends on the parameters $t_1,\dots,t_{2n+4}\in\Complex$ such that
$t_1\cdots t_{2n+4}=pq$ and
$\max\{\abs{t_1},\dots,\abs{t_{2n+4}}\}<1$, and can be stated as
\begin{multline}\label{Cnbeta}
\frac{1}{(2\pi\im)^n} \int_{\T^n}
\prod_{j=1}^n\frac{\prod_{i=1}^{2n+4}\Gamma(t_iz_j^{\pm})}
{\Gamma(z_j^{\pm 2})}
\prod_{1\leq i<j\leq n} \frac{1}{\Gamma(z_i^{\pm}z_j^{\pm})}
\d{z} \\
=\frac{2^n n!}{\pq{n}}
\prod_{1\leq i<j\leq 2n+4}\Gamma(t_it_j).
\end{multline}
This was conjectured by van Diejen and Spiridonov~\cite[Theorem 4.1]{vDS01}
who gave a proof based on a certain vanishing hypothesis and
proven in full by Rains~\cite[Corollary 3.2]{Rains03}
and the first author \cite[Theorem 2]{Spiridonov04}.

In the limit when $p$ tends to $0$ the type I A$_n$ and C$_n$
elliptic beta integrals reduce to multiple integrals
of Gustafson \cite[Theorems 2.1 and 4.1]{Gustafson92}.

The identification with the A$_n$ and C$_n$ root systems
in the above two integrals is simple.
In the case of A$_n$ the set of roots $\Delta^{\A}$ is given by
$\Delta^{\A}=\{\epsilon_i-\epsilon_j|~i,j\in\set{n+1}, i\neq j\}$ with
$\epsilon_i$ the $i$th standard unit vector in $\R^{n+1}$.
Setting $\phi_i=\epsilon_i-(\epsilon_1+\cdots+\epsilon_{n+1})/(n+1)$,
we formally put $z_i=\exp(\phi_i)$. Hence $z_1\cdots z_{n+1}=1$,
and the permutation symmetry in the $z_i$, $i\in\set{n+1}$ of the
integrand in \eqref{Anbeta}
is in accordance with the A$_n$ Weyl group, which acts on $\Delta^{\A}$ by
permuting the indices of the $\epsilon_i$.
The factor $\prod_{1\leq i<j\leq n+1}\Gamma(z_iz_j^{-1},z_i^{-1}z_j)$ in
\eqref{Anbeta} is identified with $\prod_{\alpha\in\Delta^{\A}}\Gamma(\e^{\alpha})$.

In the case of C$_n$ the set of roots is given by
$\Delta^{\C}=\{\pm 2\varepsilon_i| i\in\set{n}\}
\cup\{\pm\varepsilon_i\pm \varepsilon_j|~1\leq i<j\leq n\}$
with $\varepsilon_i$ the $i$th standard unit vector in $\R^n$ and the two $\pm$'s in $
\pm\varepsilon_i\pm \varepsilon_j$ taken independently.
Furthermore, $z_i=\exp(\varepsilon_i)$, and the hyperoctahedral (i.e., signed permutation)
symmetry of the
integrand in \eqref{Cnbeta} reflects the C$_n$ Weyl group symmetry of $\Delta^{\C}$.
The factor $\prod_{i=1}^n\Gamma(z_i^{\pm 2})\prod_{1\leq i<j\leq n}\Gamma(z_i^{\pm}z_j^{\pm})$ in
\eqref{Cnbeta} is now identified with $\prod_{\alpha\in\Delta^{\C}}\Gamma(\e^{\alpha})$.

\section{Inversion formulas. I. The single variable case}

\subsection{Motivation}\label{secmot}
To explain the origin of the inversion formula given in
Theorem~\ref{inversionthm} below let us take the elliptic beta integral
\eqref{elbetasymm} and remove the restrictions \eqref{cond}.
The price to be paid is that the contour $\T$ has to be replaced by $C$,
where $C$ is a contour\footnote{When dealing with one-dimensional contour
integrals we always assume the contour $C$ to be a positively oriented
Jordan curve such that $C=C^{-1}$, i.e., such that if $z\in C$ then also
$z^{-1}\in C$. Consequently, if a point $z$ lies in the interior of $C$ then
its reciprocal $z^{-1}$ lies in the exterior of $C$.}
such that the sequences of poles of the integrand converging to zero
(i.e., the poles at
$z=t_ip^{\mu}q^{\nu}$ for $i\in\set{6}$)
lie in the interior of $C$.
Defining
\begin{equation}\label{kappadef}
\kappa=\frac{\pqinf}{4\pi\im}
\end{equation}
we thus have
\begin{equation}\label{ebcont}
\kappa \int_C \frac{\prod_{i=1}^6\Gamma(t_iz^{\pm})}
{\Gamma(z^{\pm 2})}\d{z} \\
=\prod_{1\leq i<j\leq 6}\Gamma(t_it_j)
\end{equation}
for $t_i\in\Complex$ such that $t_1\cdots t_6=pq$.

Making the substitutions
\begin{equation}\label{sub}
(t_1,\dots,t_6)\to (t^{-1}w,t^{-1}w^{-1},ts_1,\dots,ts_4)
\end{equation}
we obtain
\begin{equation*}
\kappa \int_C \frac{\Gamma(t^{-1}w^{\pm}z^{\pm})
\prod_{i=1}^4\Gamma(ts_iz^{\pm})}
{\Gamma(z^{\pm 2})}\d{z}
=\Gamma(t^{-2})\prod_{i=1}^4\Gamma(s_iw^{\pm})
\prod_{1\leq i<j\leq 4}\Gamma(t^2s_is_j)
\end{equation*}
with $t^2s_1\cdots s_4=pq$ and $C$ a contour that has the poles
of the integrand at
\begin{equation}\label{C}
z=t^{-1}w^{\pm}p^{\mu}q^{\nu} \text{ and } z=ts_ip^{\mu}q^{\nu},\quad
i\in\set{4}
\end{equation}
in its interior.
Multiplying both sides by
$\kappa\,\Gamma(tw^{\pm}x^{\pm})/\Gamma(w^{\pm 2})$
and integrating $w$ along a contour $\hat{C}$ around
\begin{equation}\label{Chat}
w=tx^{\pm}p^{\mu}q^{\nu} \text{ and } w=s_ip^{\mu}q^{\nu},\quad
i\in\set{4}
\end{equation}
we get
\begin{align*}
\kappa^2 & \int_{\hat{C}} \int_C
\frac{\Gamma(tw^{\pm}x^{\pm},t^{-1}w^{\pm}z^{\pm})
\prod_{i=1}^4\Gamma(ts_iz^{\pm})}
{\Gamma(z^{\pm 2},w^{\pm 2})} \d{z} \d{w} \\
&=\kappa\, \Gamma(t^{-2})\prod_{1\leq i<j\leq 4}\Gamma(t^2s_is_j)
\int_{\hat{C}}
\frac{\Gamma(tw^{\pm}x^{\pm})\prod_{i=1}^4\Gamma(s_iw^{\pm})}
{\Gamma(w^{\pm 2})} \d{w} \\
&=\Gamma(t^{\pm 2})\prod_{i=1}^4\Gamma(ts_ix^{\pm}).
\end{align*}
Here the second equality follows by application of the elliptic
beta integral \eqref{ebcont}.

Inspection of the left and right-hand sides of the above result
reveals that for
\begin{equation}\label{fz}
f(z)=\prod_{i=1}^4\Gamma(ts_iz^{\pm}), \qquad t^2s_1\cdots s_4=pq
\end{equation}
the following reproducing double integral holds
\begin{equation}\label{invf}
\frac{\kappa^2}{\Gamma(t^{\pm 2})}\int_{\hat{C}} \int_C
\frac{\Gamma(tw^{\pm}x^{\pm},t^{-1}w^{\pm}z^{\pm})}
{\Gamma(z^{\pm 2},w^{\pm 2})} f(z) \d{z} \d{w} =f(x)
\end{equation}
provided the contours $C$ and $\hat{C}$ are chosen in accordance with
\eqref{C} and \eqref{Chat}.

\subsection{The $n=1$ integral inversion}\label{secn1}

If we choose $\abs{t}<1$ and $\max\{\abs{s_1},\dots,\abs{s_4}\}<1$
then the function $f$ in \eqref{fz} is
free of poles for $\abs{t}\leq\abs{z}\leq \abs{t}^{-1}$.
Moreover, if we also take $\abs{t}<\abs{x}<\abs{t}^{-1}$ then all
the points listed in \eqref{Chat} have absolute value less than one,
so that we may choose $\hat{C}$ to be the unit circle $\T$.
But assuming $w\in\T$ in \eqref{C} and further demanding that
$M<\abs{t}^2$ with $M$ defined in \eqref{Mdef}, it follows that
for the above choice of parameters all the points
listed in \eqref{C} have absolute value less than $\abs{t}$ with the
exception of $z=t^{-1}w^{\pm}$.

These considerations suggest the following generalization of
\eqref{invf} to a larger class of functions.
\begin{theorem}\label{inversionthm}
Let $p,q,t\in\Complex$ such that $M<\abs{t}^2<1$.
For fixed $w\in\T$ let $C_w$ denote a contour
inside the annulus $\An=\{z\in\Complex|~
\abs{t}-\epsilon<\abs{z}<\abs{t}^{-1}+\epsilon\}$
for infinitesimally small but positive $\epsilon$, such that
$C_w$ has the points $t^{-1}w^{\pm}$ in its interior.
Let $f(z)=f(z;t)$ be a function such that $f(z)=f(z^{-1})$ and such that
$f(z)$ is holomorphic on $\An$.
Then for $\abs{t}<\abs{x}<\abs{t}^{-1}$ there holds
\begin{equation}\label{inversion}
\kappa^2 \int_{\T}\biggl(\int_{C_w}\Delta(z,w,x;t)
f(z)\d{z}\biggr) \d{w}=f(x),
\end{equation}
where
\begin{equation}\label{Delta}
\Delta(z,w,x;t)=\frac{\Gamma(tw^{\pm}x^{\pm},t^{-1}w^{\pm}z^{\pm})}
{\Gamma(t^{\pm 2},z^{\pm 2},w^{\pm 2})}.
\end{equation}
\end{theorem}
The poles of the integrand at $z=t^{-1}w^{\pm}p^{\mu}q^{\nu}$ for $(\mu,\nu)
\neq (0,0)$ are of course also in the interior of $C_w$, but since these all have
absolute value less than $\abs{t}$
(thanks to $M<\abs{t}^2$) they do not lie in $\An$.

If one drops the condition that $f(z)=f(z^{-1})$ then
the right hand side of \eqref{inversion} should be symmetrized,
giving $(f(x)+f(x^{-1})/2$ instead of $f(x)$.

Since the kernel $\Delta(z,w,x;t)$ factorizes as
$\Delta(z,w,x;t)=\delta(z,w;t^{-1})\delta(w,x;t)$ with
\begin{equation*}
\delta(z,w;t)=\frac{\Gamma(tw^{\pm}z^{\pm})}{\Gamma(t^2,z^{\pm 2})}
\end{equation*}
the identity \eqref{inversion}
may also be put as the following elliptic integral transform. If
\begin{subequations}\label{epair}
\begin{equation}\label{fhat}
\hat{f}(w;t)=\kappa \int_{C_w} \delta(z,w;t^{-1}) f(z;t)\d{z}
\end{equation}
then
\begin{equation}\label{f}
f(x;t)=\kappa \int_{\T}\delta(w,x;t) \hat{f}(w;t)\d{w}
\end{equation}
\end{subequations}
provided all the conditions and definitions of
Theorem~\ref{inversionthm} are assumed.
The theorem may thus be formally viewed as the inversion of
the integral operator
$\delta(w;t)$ defined by
\begin{equation*}
\delta(w;t)f=\kappa
\int_{C_w} \delta(z,w;t^{-1}) f(z)\d{z}.
\end{equation*}
The external variable $w$ enters the kernel  $\delta(z,w;t)$ through the
term $\Gamma(tw^\pm z^\pm)$, which reflects only a part
of the elliptic beta integral structure \eqref{ebcont}. In this sense, we have
a universal integral transformation, playing a central role in the context
of integral Bailey pairs \cite{Spiridonov04a}. In particular, after
taking the
limit $p\to 0$, we obtain a $q$-hypergeometric integral transformation
which does not distinguish the Askey-Wilson and Rahman integrals.
In this respect, our integral transformation essentially differs from
the one introduced in \cite{KS99} on the basis of the full kernel of the
Askey-Wilson integral.

An example of a pair $(f,\hat{f})$ is given by $f$ of \eqref{fz} and
\begin{equation*}
\hat{f}(z)=\prod_{i=1}^4\Gamma(s_iz^{\pm})\prod_{1\leq i<j\leq 4}\Gamma(t^2s_is_j).
\end{equation*}
For later comparison we eliminate $s_4$ and apply \eqref{Gammasymm}. After
normalizing the above pair of functions we find the new pair
\begin{equation}\label{newpair}
f(z)=\frac{\prod_{i=1}^3\Gamma(S s_i^{-1},ts_iz^{\pm})}{\Gamma(tSz^{\pm})}
\quad\text{and}\quad
\hat{f}(z)=
\frac{\prod_{i=1}^3\Gamma(t^2Ss_i^{-1},s_iz^{\pm})}{\Gamma(t^2S z^{\pm})}
\end{equation}
with $S=s_1s_2s_3$ and $\max\{\abs{s_1},\abs{s_2},\abs{s_3},
\abs{t^{-2}S^{-1}pq}\}<1$.
Writing $f(z;t,s)$ and $\hat{f}(z;t,s)$ instead of $f(z)$ and $\hat{f}(z)$, gives
\begin{equation}\label{ffhatsymm}
\hat{f}(z;t,s)=f(z^{-1};t^{-1},ts) \quad\text{and}\quad
f(z;t,s)=\hat{f}(z^{-1};t^{-1},ts)
\end{equation}
with $s=(s_1,s_2,s_3)$ and $ts=(ts_1,ts_2,ts_3)$.
The reason for writing $z^{-1}$ and not $z$ on the
right is that it is the above form that generalizes to A$_n$, see Section~\ref{subsecIi}.
For the pair $(f,\hat{f})$ of \eqref{newpair} we can also deform the
respective contours of integration in \eqref{epair} and more
symmetrically write
\begin{equation*}
\hat{f}(z;t,s)=\kappa \int_{C_{w;t^{-1},s}} \delta(z,w;t^{-1}) f(z;t,s)\d{z}
\end{equation*}
and
\begin{equation*}
f(z;t,s)=\kappa \int_{C_{w;t,ts}}\delta(z,w;t) \hat{f}(w;t,s)\d{z},
\end{equation*}
with $C_{w;t,s}$ a contour that has the points
$t^{-1}w^{\pm}p^{\mu}q^{\nu}$, $ts_ip^{\mu}q^{\nu}$ and $t^{-1}S^{-1}p^{\mu+1}q^{\nu+1}$
in its interior, and where $t$ and $s=(s_1,s_2,s_3)$ can be chosen freely.
This can also be captured in just a single equation as
\begin{equation*}
f(z;t,s)=\kappa \int_{C_{w;t,ts}} \delta(z,w;t) f(z;t^{-1},ts)\d{z}.
\end{equation*}

\subsection{Proof of Theorem~\ref{inversionthm}}
Consider the integral over $z$ in \eqref{inversion} for
fixed $w\in\T$ such that $w^2\neq 1$.
By deforming the integration contour from $C_w$ to $\T$
the simple poles at $z=t^{-1}w^{\pm}$
($tw^{\pm}$) move from the interior (exterior) to the
exterior (interior) of the contour of integration.
Calculating the respective residues using the $f(z)=f(z^{-1})$
and $\Delta(z,w,x;t)=\Delta(z^{-1},w,x;t)$ symmetries and
the limit \eqref{Glim}, yields
\begin{multline}\label{exp1}
\kappa\int_{C_w}\Delta(z,w,x;t)f(z)\d{z}
=\kappa\int_{\T}\Delta(z,w,x;t)f(z)\d{z} \\
+\frac{\Gamma(tw^{\pm} x^{\pm})}{\Gamma(t^2)}
\biggl(\frac{f(t^{-1}w)}{\Gamma(w^2,t^2w^{-2})}
+\frac{f(t^{-1}w^{-1})}{\Gamma(w^{-2},t^2w^2)}\biggr).
\end{multline}
Since $1/\Gamma(1)=\Gamma(pq)=0$ both sides vanish identically
for $w^2=1$ so that the above is true for all $w\in\T$.

Next, by \eqref{exp1},
\begin{align*}
I(x;t)&:=\kappa^2\int_{\T}\int_{C_w}\Delta(z,w,x;t)
f(z)\d{z}\d{w} \\
&\hphantom{:}=\kappa^2\int_{\T^2}\Delta(z,w,x;t)f(z)\d{z}\d{w}
+2\kappa
\int_{\T}\frac{\Gamma(tw^{\pm} x^{\pm})}
{\Gamma(t^2,w^2,t^2w^{-2})}f(t^{-1}w) \d{w},
\end{align*}
where we have made the substitution $w\to w^{-1}$
in the integral over $w$ corresponding to the last
term on the right of \eqref{exp1}.

To proceed we replace $w\to tz$ in the single integral on the right
and invoke Fubini's theorem to interchange the order of integration
in the double integral. Hence
\begin{equation*}
I(x;t)=\kappa^2\int_{\T^2}\Delta(z,w,x;t)f(z)\d{w}\d{z}
+2\kappa \int_{t^{-1}\T}
\frac{\Gamma(x^{\pm}z^{-1},t^2x^{\pm}z)}{\Gamma(t^2,z^{-2},t^2z^2)}
f(z)\d{z}.
\end{equation*}
where $a\T$ denotes the positively oriented circle of radius $\abs{a}$.
If we deflate $t^{-1}\T$ to $\T$ the pole at $z=x$
(if $1<\abs{x}<\abs{t}^{-1}$) or $z=x^{-1}$
(if $\abs{t}<\abs{x}<1$) moves from the interior to the exterior of
the integration contour.
By the symmetry of $f$ we find
\begin{multline}\label{Iinterm}
I(x;t)=\kappa^2\int_{\T^2}\Delta(z,w,x;t)f(z)\d{w}\d{z} \\
+2\kappa\int_{\T}
\frac{\Gamma(x^{\pm}z^{-1},t^2x^{\pm}z)}{\Gamma(t^2,z^{-2},t^2z^2)}f(z)\d{z}
+f(x)
\end{multline}
irrespective of whether $\abs{t}<\abs{x}<1$ or $1<\abs{x}<\abs{t}^{-1}$.

When $\abs{x}=1$ we require the Sokhotsky--Plemelj definition of the
Cauchy integral in the case of a pole singularity
on the integration contour $C$:
\begin{equation*}
\int_C\frac{f(z)}{z-x}\,\textup{d}z=\frac{1}{2}
\int_{C_{+}}\frac{f(z)}{z-x}\,\textup{d}z+\frac{1}{2}
\int_{C_{-}}\frac{f(z)}{z-x}\,\textup{d}z,
\end{equation*}
where $f(z)$ is holomorphic on $C$, and $C_{\pm}$ are contours
which include/exclude the point $x\in C$ by an
infinitesimally small deformations of $C$ in the vicinity of $x$.
By the $x\to x^{-1}$ symmetry of our integral it thus follows that
\eqref{Iinterm} is true for all $\abs{t}<\abs{x}<\abs{t}^{-1}$.

To complete the proof we need to show that the integrals
on the right-hand side of \eqref{Iinterm} vanish.
To achieve this we use that for $z\in\T$
\begin{multline}\label{exp2}
\kappa\int_{\T}\Delta(z,w,x;t)\d{w} \\
=\kappa\int_C\Delta(z,w,x;t)\d{w}
-\frac{\Gamma(x^\pm z^{-1},t^2x^{\pm}z)}{\Gamma(t^2,z^{-2},t^2z^{2})}
-\frac{\Gamma(x^\pm z,t^2x^{\pm} z^{-1})}{\Gamma(t^2,z^2,t^2z^{-2})},
\end{multline}
where $C$ is a contour such that the points $w=tx^{\pm} p^{\mu}q^{\nu}$
and $w=t^{-1}z^{\pm} p^{\mu}q^{\nu}$ lie in its interior.
The two ratios of elliptic gamma functions on the right correspond to
the residues of the poles at $w=t^{-1}z^{\pm}$ and $w=tz^{\pm}$
which, for $\abs{z}=1$ and $\abs{t}<1$, lie in the exterior and interior
of $\T$, respectively.
Note that we again have implicitly assumed $z^2\neq 1$ in the calculation
of the respective residues, but that \eqref{exp2} is true for all $z\in\T$.

Since $\Gamma(pq)=0$ it follows from the elliptic beta
integral \eqref{ebcont} with $t_5t_6=pq$ that the integral
on the right vanishes, resulting in
\begin{equation*}
\kappa\int_{\T}\Delta(z,w,x;t)\d{w}
=-\frac{\Gamma(x^\pm z^{-1},t^2x^{\pm}z)}{\Gamma(t^2,z^{-2},t^2z^{2})}
-\frac{\Gamma(x^\pm z,t^2x^{\pm} z^{-1})}{\Gamma(t^2,z^2,t^2z^{-2})}.
\end{equation*}
Substituting this in the first term on the right of \eqref{Iinterm}
and making a $z\to z^{-1}$ variable change establishes the desired cancellation
of integrals in \eqref{Iinterm}, thereby establishing the theorem.

\section{Inversion formulas. II. The root systems A$_n$ and C$_n$}

\subsection{Main results}\label{secmain}
To state our multi-dimensional inversion theorems
we first extend the (what will be referred to as A$_1$ or C$_1$)
symmetry $f(z)=f(z^{-1})$ to functions of $n$ variables.
Let $g$ be a symmetric function of $n+1$ independent variables.
Then a function $f(z)=f(z_1,\dots,z_n):=g(z_1,\dots,z_{n+1})$ is said to have
A$_n$ symmetry. (Recall our convention that $z_1\cdots z_{n+1}=1$.)
Similarly, we say that $f(z)=f(z_1,\dots,z_n)$ has C$_n$
symmetry if $f$ is symmetric under signed permutations.
That is, $f(z)=f(\varw(z))$ for $\varw\in S_n$ and
$f(z_1,\dots,z_n)=f(z_1^{\sigma_1},\dots,z_n^{\sigma_n})$
where each $\sigma_i\in\{-1,1\}$.
For example $f$ has A$_2$ symmetry if
$f(z_1,z_2)=f(z_2,z_1)$ and
$f(z_1,z_2)=f(z_1,z_3)=f(z_1,z_1^{-1}z_2^{-1})$,
and $f$ has C$_2$ symmetry if
$f(z_1,z_2)=f(z_2,z_1)$ and $f(z_1,z_2)=f(z_1,z_2^{-1})$.
The integrands of the integrals \eqref{Anbeta} and \eqref{Cnbeta}
provide examples of functions that are A$_n$ or C$_n$ symmetric.

Below we will also use the root system analogues of $\kappa$
of equation~\eqref{kappadef};
\begin{equation}\label{kappan}
\kappa^{\A}=\frac{\pq{n}}{(2\pi\im)^n(n+1)!}
\quad\text{and}\quad
\kappa^{\C}=\frac{\pq{n}}{(2\pi\im)^n2^n n!}.
\end{equation}

Finally, we need to discuss a somewhat technically involved issue.
The $n=1$ inversion formula \eqref{inversion} features the
integration contour $C_w$ which is a deformation of the
contour $\T$ such that the poles of the integrand at
$t^{-1}w^{\pm}p^{\mu}q^{\nu}$ are in the interior of $C_w$.
Now the A$_n$ beta integral \eqref{Anbeta} is computed by iteratively
integrating over the $n$ components of $z$.
Let us choose to integrate $z_n$ first then $z_{n-1}$ and so on.
When doing the $z_i$ integral the integrand will have poles
which are independent of $z_1,\dots,z_{i-1}$ and poles which
depend on these variables through their product
$Z_{i-1}:=z_1\cdots z_{i-1}$.
For example, when doing the $z_n$ integral over $\T$ we need
to compute the residues of the poles at $z_n=t_ip^{\mu}q^{\nu}$
and $z_{n+1}=s_i^{-1}p^{-\mu}q^{-\nu}$, i.e., at
$z_n=s_ip^{\mu}q^{\nu}Z_{n-1}^{-1}$.
Just as in the $n=1$ case we wish to utilize the A$_n$ beta
integral in which $(t_1,\dots,t_{n+1})$ is substituted by
$(t^{-1}w_1,\dots,t^{-1}w_{n+1})$ with $w_i\in\T$. (Compare this
with \eqref{sub}.)
Hence we need to again analytically continue
the integral \eqref{Anbeta} by appropriately deforming the
integration contours. Because of the above-discussed poles
depending on the remaining integration variables this
deformation --- which will be denoted by $C_w^n$ ---
cannot be of the form $C_1\times C_2\times\dots
\times C_n$ with each of the one-dimensional contours
$C_i$ independent of $z$.
Rather what we get is that $C_n$ depends on $t$ and $w$ as well as on
$Z_{n-1}$. Then $C_{n-1}$ will depend on $t$, $w$, and $Z_{n-2}$ and
so on.
Of course, this is all assuming the above order of integrating out the
components of $z$, but, evidently, all ordering are in fact equivalent.

We would like an efficient description of the deformed contours
that is independent of the chosen order of integration
and that reflects the A$_n$ symmetry present in the problem.
However, since we want to avoid the complexities of genuine higher-dimensional
residue calculus, we adopt a convention that $C_w^n$ does not
explicitly describe each of the one-dimensional contours
composing it. Rather, we encode $C_w^n$ by indicating which
poles of the integrand are to be taken in the interior and exterior
at each stage of the iterative computation of the integral over $z$.

Let $p,q,t\in\Complex$ such that $M<\abs{t}^{n+1}<1$ and denote
\begin{equation}\label{An}
\An=\{z\in\Complex^n|~\abs{z_j}<\abs{t}^{-1}+\epsilon,~
j\in\set{n+1}\}
\end{equation}
for infinitesimally small but positive $\epsilon$.
Let $f$ be an A$_n$ symmetric function holomorphic on $\An$
and let the generalization of the kernel \eqref{Delta}
to the root system pair (A$_n$,A$_n$) be given by
\begin{multline}\label{DeltaA}
\Delta^{(\A,\A)}(z,w,x;t) \\
=\frac{\prod_{i,j=1}^{n+1} \Gamma(tw_ix_j,t^{-1}w_i^{-1}z_j^{-1})}
{\Gamma(t^{n+1},t^{-n-1})
\prod_{1\leq i<j\leq n+1}
\Gamma(z_iz_j^{-1},z_i^{-1}z_j,w_iw_j^{-1},w_i^{-1}w_j)}.
\end{multline}
Then for $w\in\T^n$ we write
\begin{equation}\label{Cint}
\int_{C^n_w} \Delta^{(\A,\A)}(z,w,x;t)f(z)\d{z},
\end{equation}
where --- by abuse of notation --- we write
`$C^n_w\subset\An$' as the `deformation of the
\textup{(}oriented\textup{)} $n$-torus $\T^n$' such that for all
$i\in\set{n+1}$
\begin{equation}\label{poles}
z_j=t^{-1}w_i^{-1}~
\begin{cases}
\text{lies in the interior of $C^n_w$ if $j\in\set{n}$} \\
\text{lies in the exterior of $C^n_w$ if $j=n+1$}.
\end{cases}
\end{equation}

More precisely, we consider $C_w^n$ as an iteratively
defined $n$-dimensional structure encoding which poles of
the integrand are to be taken in the
interior/exterior at each stage of the iterative integration over $z$.
That is, if we again fix the order of integration as before then,
when integrating over $z_j$, the poles (these will occur regardless
of which components of $z$ are already integrated out) at
$z_j=t^{-1}w_i^{-1}p^{\mu}q^{\nu}$ are all in the interior of $C_j$
because (i) for $(\mu,\nu)\neq (0,0)$ we have
$\abs{t^{-1}w_i^{-1}p^{\mu}q^{\nu}}<\abs{t}^{-1}M<\abs{t}^n$,
but for $z\in\An$ each $z_j$ is bounded (in absolute value)
from below by $\abs{t}^n$,
(ii) for $(\mu,\nu)=(0,0)$ we have to satisfy \eqref{poles}.
The poles at $z_j=t^{n-j+1}w_{i_1}\cdots w_{i_{n-j+1}}p^{-\mu}q^{-\nu}Z_{j-1}^{-1}$
(corresponding to the pole $z_{n+1}=t^{-1}w_i^{-1}p^{\mu}q^{\nu}$
with $z_n,\dots,z_{j+1}$ integrated out)
will be in the exterior of $C_j$ because
(i) for $(\mu,\nu)\neq (0,0)$ we have
$\abs{t^{n-j+1}w_{i_1}\cdots w_{i_{n-j+1}}p^{-\mu}q^{-\nu}Z_{j-1}^{-1}}>
\abs{t^{n-j+1}Z_{j-1}^{-1}}M^{-1}>\abs{t^{-j}Z_{j-1}^{-1}}>\abs{t}^{-1}$
(by \eqref{An}) and $C_w^n\subset\An$),
but for $z\in\An$, $\abs{z_j}<\abs{t}^{-1}$,
(ii) for $(\mu,\nu)=(0,0)$ we have to satisfy \eqref{poles}.
Of course, because $f$ is holomorphic on $\An$ its poles are either
trivially in the interior or exterior of each contour $C_j$.

The most rigorous definition of the integration domain $C_w^n\subset\An$ is
obtained by considering it as a deformation of $\T^n$ allowing for an
analytical  continuation of the integral
\begin{equation*}
\int_{\T^n} \frac{\prod_{i,j=1}^{n+1}\Gamma(t_iz_j^{-1})}
{\prod_{1\leq i<j\leq n+1} \Gamma(z_iz_j^{-1},z_i^{-1}z_j)}\: f(z)\frac{dz}{z}
\end{equation*}
from the restricted values of parameters $|t_i|<1$, $i\in [n+1]$,
to the region $t_i=t^{-1}w_i^{-1}$ with $|t|<1$ and $|w_i|=1,$
$w_1\cdots w_{n+1}=1$. Clearly, this defines the
$z$-dependent part of the integral \eqref{Cint}.
 However, for making our computations efficient we will not
 reformulate our results in this coordinate independent way
 but characterize $C_w^n$ by locations of the appropriate poles.

With a trivial modification of the above notation we can now
formulate two generalizations of \eqref{inversion} corresponding
to the root system pairs
$(\text{A}_n,\text{A}_n)$ and $(\text{A}_n,\text{C}_n)$.
In the next section we shall
also formulate a conjecture for the pair $(\text{C}_n,\text{A}_n)$.

\begin{theorem}[\textbf{(A,A) inversion formula}]\label{inversionthmA}
Let $q,p,t\in\Complex$ such that $M<\abs{t}^{n+1}<1$ and let
$\An$ be defined as in \eqref{An}.
For fixed $w\in\T^n$ let $C^n_w\subset\An$ denote a deformation of the
\textup{(}oriented\textup{)} $n$-torus $\T^n$ such that
\eqref{poles} holds for all $i\in\set{n+1}$.
Let $f$ be an A$_n$ symmetric function holomorphic on $\An$.
Then for $x\in\Complex^n$ such that
$\abs{x_j}<\abs{t}^{-1}$ for all $j\in\set{n+1}$ and
\begin{equation}\label{extra}
\abs{x_j}>1 \quad \text{for $j\in\set{n}$},
\end{equation}
there holds
\begin{equation}\label{invA}
(\kappa^{\A})^2
\int_{\T^n}\biggl(\int_{C^n_w}
\Delta^{(\A,\A)}(z,w,x;t)f(z)\d{z}\biggr)\d{w}=f(x),
\end{equation}
where $\Delta^{(\A,\A)}(z,w,x;t)$ is given in \eqref{DeltaA}.
\end{theorem}
Because of the A$_n$ symmetry in $x$ of \eqref{invA}
the condition \eqref{extra} can of course be replaced by
the condition that all but one of $\abs{x_1},\dots,\abs{x_{n+1}}$
exceeds one.
In fact, we have strong evidence that the condition \eqref{extra} is
not necessary. However, the proof of the theorem becomes significantly
more complicated if \eqref{extra} is dropped and in the absence of
\eqref{extra} we have only been able to complete the proof for $n\leq 2$.

\begin{theorem}[\textbf{(A,C) inversion formula}]\label{inversionthmAC}
Let $p,q,t\in\Complex$ such that $M<\abs{t}^{n+1}<1$ and let
$\An$ be defined as in \eqref{An}.
Let $C_n^w\subset\An$ be a deformation of $\T^n$ such that
for all $i\in\set{n}$
\begin{equation}\label{poles2}
z_j=t^{-1}w_i^{\pm}~
\begin{cases}
\text{lies in the interior of $C^n_w$ if $j\in\set{n}$} \\
\text{lies in the exterior of $C^n_w$ if $j=n+1$}.
\end{cases}
\end{equation}
Let $f$ be an A$_n$ symmetric function holomorphic on $\An$.
Then for $x\in\Complex^n$ such that
$\abs{x_j}<\abs{t}^{-1}$ for all $j\in\set{n+1}$ and
such that \eqref{extra} holds, we have
\begin{equation}\label{invAC}
\kappa^{\A}\kappa^{\C}
\int_{\T^n}\biggl(\int_{C^n_w}
\Delta^{(\A,\C)}(z,w,x;t)f(z)\d{z}\biggr)\d{w}=f(x).
\end{equation}
where
\begin{multline}\label{DeltanAC}
\Delta^{(\A,\C)}(z,w,x;t)=
\frac{\prod_{i=1}^n\prod_{j=1}^{n+1}
\Gamma(tw_i^{\pm}x_j,t^{-1}w_i^{\pm}z_j^{-1})}
{\prod_{i=1}^n\Gamma(w_i^{\pm 2})
\prod_{1\leq i<j\leq n}\Gamma(w_i^{\pm}w_j^{\pm})} \\
\times\frac{1}{\prod_{1\leq i<j\leq n+1}
\Gamma(z_iz_j^{-1},z_i^{-1}z_j,t^2z_iz_j,
t^{-2}z_i^{-1}z_j^{-1})}.
\end{multline}
\end{theorem}
Again the condition \eqref{extra} is probably unnecessary, but
without it our proof of Theorem~\ref{inversionthmAC}
requires some rather intricate modifications for $n\geq 3$.

If one drops the condition that $f$ is an
A$_n$ symmetric function then it is immediate from the
A$_n$ symmetry of the left-hand side that
$f(x)$ on the right of \eqref{invA} and \eqref{invAC}
should be replaced by the A$_n$ symmetric
\begin{equation*}
\frac{1}{(n+1)!}\sum_{\varw\in S_{n+1}}
g(x_{\varw_1},\dots,x_{\varw_{n+1}}),
\end{equation*}
where $g(x_1,\dots,x_{n+1})=f(x_1,\dots,x_n)$.

The proofs of the two inversion theorems are very similar,
the only significant difference being that \eqref{invA} requires the
A$_n$ elliptic beta integral and \eqref{invAC} the C$_n$ elliptic
beta integral to establish the vanishing of certain unwanted
terms arising in the expansion of the integral over $C^n_w$ as
a sum of integrals over $\T^n,\dots,\T,\T^0$.
We therefore content ourselves with only presenting the details of the proof of
Theorem~\ref{inversionthmAC}.

First, however, let us state the root systems analogues
of some of the equations of Section~\ref{secn1}.
This will lead us to discover a new elliptic beta integral
for the root system A$_n$. Loosely speaking this new integral may be
viewed as the inverse of the C$_n$ beta integral \eqref{Cnbeta}
with respect to the kernel $\Delta^{(\A,\C)}$.

\subsection{Consequences of Theorems~\ref{inversionthmA} and
\ref{inversionthmAC}}\label{seccons}

\subsubsection{Integral identities}\label{subsecIi}

Introducing
\begin{equation}\label{deltaA}
\nabla^{\A}(z,w;t)=
\frac{\prod_{i,j=1}^{n+1} \Gamma(tw_i^{-1}z_j^{-1})}
{\Gamma(t^{n+1})
\prod_{1\leq i<j\leq n+1}
\Gamma(z_iz_j^{-1},z_i^{-1}z_j)},
\end{equation}
the kernel $\Delta^{(\A,\A)}(z,w,x;t)$ may be factored as
\begin{equation*}
\Delta^{(\A,\A)}(z,w,x;t)=\nabla^{\A}(z,w;t^{-1})\nabla^{\A}(w^{-1},x^{-1};t),
\end{equation*}
with $w^{-1}=(w_1^{-1},\dots,w_n^{-1})$ and
$x^{-1}=(x_1^{-1},\dots,x_n^{-1})$.

Defining the elliptic integral transform
\begin{equation*}
\hat{f}^{\A}(w;t)=
\kappa^{\A} \int_{C^n_w} \nabla^{\A}(z,w;t^{-1}) f^{\A}(z;t)\d{z}
\end{equation*}
the claim of Theorem~\ref{inversionthmA} is equivalent to the
inverse transformation
\begin{equation*}
f^{\A}(x;t)=\kappa^{\A} \int_{\T^n}\nabla^{\A}(w^{-1},x^{-1};t)
\hat{f}^{\A}(w;t)\d{w}.
\end{equation*}

{}From the A$_n$ elliptic beta integral one readily obtains the following
example of a pair $(f^{\A},\hat{f}^{\A})$:
\begin{align*}
f^{\A}(z)&=\prod_{i=1}^{n+2}\Gamma(Ss_i^{-1})
\prod_{j=1}^{n+1}
\frac{\prod_{i=1}^{n+2}\Gamma(ts_iz_j)}{\Gamma(tSz_j)} \\
\hat{f}^{\A}(z)&=
\prod_{i=1}^{n+2}\Gamma(t^{n+1}Ss_i^{-1})
\prod_{j=1}^{n+1}\frac{\prod_{i=1}^{n+2}\Gamma(s_iz_j^{-1})}
{\Gamma(t^{n+1}Sz_j^{-1})},
\end{align*}
with $S=s_1\cdots s_{n+2}$ and $\max\{\abs{s_1},\dots,\abs{s_{n+2}},
\abs{t^{-n-1}S^{-1}pq}\}<1$. Here the conditions on the $s_i$ ensure that
$f^{\A}(z)$ is holomorphic on $\An$ as follows from a reasoning
similar to the one presented immediately after Theorem~\ref{inversionthm}.
We also remark that $f^{\A}$ and $\hat{f}^{\A}$ are again related by the
simple symmetry \eqref{ffhatsymm} provided we now take
$s=(s_1,\dots,s_{n+2})$.

Next we turn our attention to Theorem~\ref{inversionthmAC} and define
\begin{equation*}
\delta^{\C}(w,x;t)=
\frac{\prod_{i=1}^n\prod_{j=1}^{n+1}
\Gamma(tw_i^{\pm}x_j)}
{\prod_{i=1}^n\Gamma(w_i^{\pm 2})
\prod_{1\leq i<j\leq n}\Gamma(w_i^{\pm}w_j^{\pm})}
\end{equation*}
and
\begin{equation*}
\delta^{\A}(z,w;t)=
\frac{\prod_{i=1}^n\prod_{j=1}^{n+1}
\Gamma(tw_i^{\pm}z_j^{-1})}
{\prod_{1\leq i<j\leq n+1}
\Gamma(z_iz_j^{-1},z_i^{-1}z_j,t^{-2}z_iz_j,t^2z_i^{-1}z_j^{-1})},
\end{equation*}
so that
\begin{equation}\label{Deltadelta2}
\Delta^{(\A,\C)}(z,w,x;t)=\delta^{\A}(z,w;t^{-1})
\delta^{\C}(w,x;t).
\end{equation}
Then, according to Theorem~\ref{inversionthmAC}, if
\begin{subequations}\label{ff}
\begin{equation}\label{fwfx}
\hat{f}^{\C}(w;t)=
\kappa^{\A}
\int_{C_w^n} \delta^{\A}(z,w;t^{-1}) f^{\A}(z;t)\d{z}
\end{equation}
then
\begin{equation}\label{fxfw}
f^{\A}(x;t)=\kappa^{\C}
\int_{\T^n}\delta^{\C}(w,x;t) \hat{f}^{\C}(w;t)\d{w}.
\end{equation}
\end{subequations}
We note that for $n=1$ this simplifies to \eqref{fhat} and \eqref{f}
up to factors of $\Gamma(t^{\pm 2})$.

Let us now choose $\hat{f}^{\C}(w;t)$ as
\begin{subequations}\label{pairAC}
\begin{equation}
\hat{f}^{\C}(w;t)=\prod_{i=1}^n\prod_{j=1}^{n+3}\Gamma(w_i^{\pm}s_j)
\end{equation}
with $\max\{\abs{s_1},\dots,\abs{s_{n+3}}\}<1$
and $t^{n+1}s_1\cdots s_{n+3}=pq$.
Provided that
\begin{equation*}
\max\{\abs{x_1},\dots,\abs{x_{n+1}}\}<\abs{t}^{-1}
\end{equation*}
we can evaluate the integral \eqref{fxfw} by the C$_n$ elliptic beta integral
\eqref{Cnbeta} with $t_j\to tx_j$ for $j\in\set{n+1}$ and
$t_{j+n+1}\to s_j$ for $j\in\set{n+3}$, to find
\begin{equation}
f^{\A}(x;t)=\prod_{i=1}^{n+1}\prod_{j=1}^{n+3}\Gamma(tx_is_j)
\prod_{1\leq i<j\leq n+1} \Gamma(t^2x_ix_j)
\prod_{1\leq i<j\leq n+3} \Gamma(s_is_j).
\end{equation}
\end{subequations}
Substituting this in \eqref{fwfx} and observing that $f^{\A}$ satisfies
the conditions imposed by the theorem, we obtain the new elliptic beta integral
\begin{multline*}
\kappa^{\A} \int_{C_w^n}
\prod_{j=1}^{n+1}\frac{\prod_{i=1}^n \Gamma(t^{-1}w_i^{\pm}z_j^{-1})
\prod_{i=1}^{n+3}\Gamma(ts_iz_j)}
{\prod_{1\leq i<j\leq n+1}
\Gamma(z_iz_j^{-1},z_i^{-1}z_j,t^{-2}z_i^{-1}z_j^{-1})} \d{z} \\
=\prod_{i=1}^n\prod_{j=1}^{n+3}\Gamma(s_jw_i^{\pm})
\prod_{1\leq i<j\leq n+3} \frac{1}{\Gamma(s_is_j)}.
\end{multline*}

By appropriately deforming the contour of integration
this may be analytically continued to $\abs{t}>1$. Then replacing
$t\to t^{-1}$ and $w_i\to t_i$ and $z_i\to z_i^{-1}$,
the result can be written as an integral over $\T^n$;
\begin{multline*}
\kappa^{\A} \int_{\T^n}
\prod_{j=1}^{n+1}\frac{\prod_{i=1}^n \Gamma(tt_i^{\pm}z_j)
\prod_{i=1}^{n+3}\Gamma(t^{-1}s_iz_j^{-1})}
{\prod_{1\leq i<j\leq n+1}
\Gamma(z_iz_j^{-1},z_i^{-1}z_j,t^2z_iz_j)} \d{z} \\
=\prod_{i=1}^n\prod_{j=1}^{n+3}\Gamma(t_i^{\pm}s_j)
\prod_{1\leq i<j\leq n+3} \frac{1}{\Gamma(s_is_j)}
\end{multline*}
for $s_1\cdots s_{n+3}=t^{n+1}pq$ and $\max\{\abs{t},\abs{tt_1^{\pm}},
\dots,\abs{tt_n^{\pm}},\abs{t^{-1}s_1},\dots,\abs{t^{-1}s_{n+3}}\}<1$.
Alternatively, we may choose to replace $t_i\to t^{-1}t_i$ and
$s_i\to ts_i$ and then eliminate $t^2$ using
$t^2=pqS^{-1}$ with $S=s_1\cdots s_{n+3}$.
By \eqref{Gammasymm} this results in our next theorem.
\begin{theorem}\label{thmbetaA}
Let $t_1,\dots,t_n,s_1,\dots,s_{n+3}\in\Complex$ and
$S=s_1\cdots s_{n+3}$. For
\begin{equation}\label{stcond}
\max\{\abs{t_1},\dots,\abs{t_n},\abs{s_1},\dots,\abs{s_{n+3}},
\abs{pqS^{-1}t_1^{-1}},\dots,\abs{pqS^{-1}t_n^{-1}}\}<1
\end{equation}
there holds
\begin{multline*}
\frac{1}{(2\pi\im)^n} \int_{\T^n}
\prod_{j=1}^{n+1}
\frac{\prod_{i=1}^n \Gamma(t_iz_j)\prod_{i=1}^{n+3}\Gamma(s_iz_j^{-1})}
{\prod_{i=1}^n \Gamma(St_iz_j^{-1})}
\prod_{1\leq i<j\leq n+1}
\frac{\Gamma(Sz_i^{-1}z_j^{-1})}
{\Gamma(z_iz_j^{-1},z_i^{-1}z_j)} \d{z} \\
=\frac{(n+1)!}{\pq{n}}
\prod_{i=1}^n\prod_{j=1}^{n+3}
\frac{\Gamma(t_is_j)}{\Gamma(St_is_j^{-1})}
\prod_{1\leq i<j\leq n+3} \Gamma(Ss_i^{-1}s_j^{-1}).
\end{multline*}
\end{theorem}
We will give an independent proof of this new A$_n$ elliptic beta integral
in Section~\ref{secaltproof}.
Somewhat surprising, even when $p$ tends to zero the corresponding
beta integral is new;
\begin{multline*}
\frac{1}{(2\pi\im)^n} \int_{\T^n}
\prod_{j=1}^{n+1}
\frac{\prod_{i=1}^n (St_iz_j^{-1};q)_{\infty}}
{\prod_{i=1}^n (t_iz_j;q)_{\infty}\prod_{i=1}^{n+3}(s_iz_j^{-1};q)_{\infty}}
\prod_{1\leq i<j\leq n+1}
\frac{(z_iz_j^{-1},z_i^{-1}z_j;q)_{\infty}}
{(Sz_i^{-1}z_j^{-1};q)_{\infty}} \d{z} \\
=\frac{(n+1)!}{(q;q)_{\infty}^n}
\prod_{i=1}^n\prod_{j=1}^{n+3}
\frac{(St_is_j^{-1};q)_{\infty}}{(t_is_j;q)_{\infty}}
\prod_{1\leq i<j\leq n+3} \frac{1}{(Ss_i^{-1}s_j^{-1};q)_{\infty}}
\end{multline*}
for $\max\{\abs{t_1},\dots,\abs{t_n},\abs{s_1},\dots,\abs{s_{n+3}}\}<1$.
When $s_{n+3}$ tends to zero this reduces to a limiting case of
Gustafson's SU$(n)$ $q$-beta integral \cite[Theorem 2.1]{Gustafson92}.

Formally the further limit $q\to 1^{-}$ can be taken by
replacing $z_j\to q^{u_j}$, $t_j\to q^{a_j}$ and $s_j\to q^{b_j}$
and choosing $q=\exp(-\pi/\alpha)$ for $\alpha$ positive and real.
The integral can then be conveniently expressed in terms of the
$q$-Gamma function
\begin{equation*}
\Gamma_q(x)=(1-q)^{1-x}\frac{(q;q)_{\infty}}{(q^x;q)_{\infty}}.
\end{equation*}
By
\begin{equation*}
\Gamma(x)=\lim_{q\to 1^{-}}\Gamma_q(x),
\end{equation*}
with $\Gamma(x)$ the classical gamma function,
the limit when $\alpha$ tends to infinity is readily obtained.
\begin{theorem}\label{q=1}
Let $a_1,\dots,a_n,b_1,\dots,b_{n+3}\in\Complex$ and
$B=b_1+\dots+b_{n+3}$ such that
\begin{equation*}
\min\{\Re(a_1),\dots,\Re(a_n),\Re(b_1),\dots,\Re(b_{n+3})\}>0.
\end{equation*}
Denote by $\Gamma(x)$ the classical instead of elliptic gamma
function. Then
\begin{multline*}
\frac{1}{(2\pi\im)^n}
\int_{-\im\infty}^{\im\infty} \dots
\int_{-\im\infty}^{\im\infty}
\prod_{j=1}^{n+1}
\frac{\prod_{i=1}^n \Gamma(a_i+u_j)
\prod_{i=1}^{n+3}\Gamma(b_i-u_j)}
{\prod_{i=1}^n \Gamma(B+a_i-u_j)} \\
\qquad \times
\prod_{1\leq i<j\leq n+1}
\frac{\Gamma(B-u_i-u_j)}
{\Gamma(u_i-u_j,u_j-u_i)}\, \textup{d}u_1\cdots \textup{d}u_n \\
=(n+1)!
\prod_{i=1}^n\prod_{j=1}^{n+3}
\frac{\Gamma(a_i+b_j)}{\Gamma(B+a_i-b_j)}
\prod_{1\leq i<j\leq n+3} \Gamma(B-b_i-b_j)
\end{multline*}
with $u_1+\cdots+u_{n+1}=0$.
\end{theorem}
In the large $b_{n+3}$ limit this coincides with
the large $\alpha_n$ limit of \cite[Theorem 5.1]{Gustafson92}.
A more rigorous justification of Theorem~\ref{q=1} can be given
using the technique of Section~\ref{secaltproof}.

The above discussion of the new elliptic beta integral
suggests that Theorem~\ref{inversionthmAC} should have
the following companion.
\begin{conjecture}[\textbf{(C,A) inversion formula}]\label{conj}
Let $p,q,t\in\Complex$ such that $M<\abs{t}^{n+1}<1$.
For fixed $w\in\T^n$ let $C_w$ denote a contour inside the
annulus $\An=\{z\in\Complex|~
\abs{t}-\epsilon<\abs{z}<\abs{t}^{-1}+\epsilon\}$
for infinitesimally small but positive $\epsilon$, such that
$C_w$ has the points $t^{-1}w_j$ for $j\in\set{n+1}$ in its interior,
and set $C_w^n=C_w\times\cdots\times C_w$.
For $f$ a C$_n$ symmetric function holomorphic on $\An^n$,
and $x\in\Complex^n$ such that $\abs{t}<\abs{x_j}<\abs{t}^{-1}$,
we have
\begin{equation}\label{invCA}
\kappa^{\A}\kappa^{\C}
\int_{\T^n}\biggl(\int_{C^n_w}
\Delta^{(\C,\A)}(z,w,x;t)f(z)\d{z}\biggr)\d{w}=f(x),
\end{equation}
where
\begin{multline*}
\Delta^{(\C,\A)}(z,w,x;t)=
\frac{\prod_{i=1}^n\prod_{j=1}^{n+1}
\Gamma(tx_i^{\pm}w_j^{-1},t^{-1}z_i^{\pm}w_j)}
{\prod_{i=1}^n\Gamma(z_i^{\pm 2})
\prod_{1\leq i<j\leq n}\Gamma(z_i^{\pm}z_j^{\pm})} \\
\times\frac{1}{
\prod_{1\leq i<j\leq n+1}\Gamma(w_iw_j^{-1},w_i^{-1}w_j,
t^{-2}w_iw_j,t^2w_i^{-1}w_j^{-1})}.
\end{multline*}
\end{conjecture}
Note that
\begin{equation*}
\Delta^{(\C,\A)}(z,w,x;t)=\delta^{\C}(z,w;t^{-1})
\delta^{\A}(w,x;t)
\end{equation*}
to be compared with \eqref{Deltadelta2}.
Hence, if this conjecture were true then
\begin{subequations}
\begin{equation}\label{ff1}
\hat{f}^{\A}(w;t)=
\kappa^{\C}
\int_{C_w^n}\delta^{\C}(z,w;t^{-1}) f^{\C}(z;t)\d{z}
\end{equation}
would imply
\begin{equation}\label{ff2}
f^{\C}(x;t)=\kappa^{\A}
\int_{\T^n}\delta^{\A}(w,x;t) \hat{f}^{\A}(w;t)\d{w}
\end{equation}
\end{subequations}
which is the image of \eqref{ff} under the interchange of the
root systems A and C.
If we take
\begin{subequations}\label{pairCA}
\begin{equation}
f^{\C}(z;t)=\prod_{i=1}^n\prod_{j=1}^{n+3}\Gamma(ts_jz_i^{\pm})
\end{equation}
for $\max\{\abs{s_1},\dots,\abs{s_{n+3}}\}<1$ and $t^2s_1\cdots s_{n+3}=pq$,
then the integral \eqref{ff1} can be calculated using a deformation
of the C$_n$ elliptic beta integral \eqref{Cnbeta}, and
\begin{equation}
\hat{f}^{\A}(w;t)=\prod_{i=1}^{n+1}\prod_{j=1}^{n+3}\Gamma(w_is_j)
\prod_{1\leq i<j\leq n+1} \Gamma(t^{-2}w_iw_j)
\prod_{1\leq i<j\leq n+3} \Gamma(t^2s_is_j).
\end{equation}
\end{subequations}
Substituting this in \eqref{ff2} once more yields the integral
of Theorem~\ref{thmbetaA} (up to simple changes of variables).

Eliminating $s_{n+3}$ in the functions listed in \eqref{pairAC} and
\eqref{pairCA}, and making the $s$ dependence explicit, we get
the formal symmetry relations
\begin{equation*}
\hat{f}^{\C}(z;t,s)=f^{\C}(z;t^{-1},ts)
\quad \text{ and } \quad
\hat{f}^{\A}(z;t,s)=f^{\A}(z;t^{-1},ts)
\end{equation*}
where $s=(s_1,\dots,s_{n+2})$ and $ts=(ts_1,\dots,ts_{n+2})$.

As already mentioned at the end of Section~\ref{secmain},
the proofs of Theorems~\ref{inversionthmA} and
\ref{inversionthmAC} hinge on the vanishing of certain
unwanted elliptic hypergeometric integrals. Key to this are specializations
of the the A$_n$ and C$_n$ elliptic beta integrals.
If a similar approach is taken with respect to
Conjecture~\ref{conj} one encounters unwanted integrals
for which the vanishing condition is not evident.

\subsubsection{Series identities}

Using residue calculus one can reduce elliptic beta integrals to
summation identities for elliptic hypergeometric series, see e.g.,
\cite{vDS00,Spiridonov03b} for examples of this procedure.
Below we will give the main steps of such a calculation
for the elliptic beta integral of Theorem~\ref{thmbetaA}.

First let us write the integral in question in the form
\begin{equation}\label{rhoint}
\kappa^{\A}\int_{\T^n} \rho(z;s,t) \d{z}=1,
\end{equation}
with integration kernel $\rho(z;s,t)$ for $s\in\Complex^{n+3}$
and $t\in\Complex^n$ given by
\begin{multline}\label{rhodef}
\rho(z;s,t)=
\prod_{j=1}^{n+1}
\frac{\prod_{i=1}^n \Gamma(t_iz_j)\prod_{i=1}^{n+3}\Gamma(s_iz_j^{-1})}
{\prod_{i=1}^n \Gamma(St_iz_j^{-1})}
\prod_{1\leq i<j\leq n+1}
\frac{\Gamma(Sz_i^{-1}z_j^{-1})}
{\Gamma(z_iz_j^{-1},z_i^{-1}z_j)} \\
\times
\prod_{i=1}^n\prod_{j=1}^{n+3}
\frac{\Gamma(St_is_j^{-1})}{\Gamma(t_is_j)}
\prod_{1\leq i<j\leq n+3}\frac{1}{\Gamma(Ss_i^{-1}s_j^{-1})}.
\end{multline}

The following poles of $\rho(z;s,t)$ lie in the interior of $\T^n$:
\begin{align}\label{inset}
z_{n+1}^{-1}&=t_ip^{\mu}q^{\nu} && i\in\set{n} \notag \\
z_j&=s_ip^{\mu}q^{\nu} && i\in\set{n+3},~j\in\set{n} \\
z_{n+1}^{-1}&=S^{-1}t_i^{-1}p^{\mu+1}q^{\nu+1}
\hspace{-15mm} && i\in\set{n} \notag \\
z_iz_j&=Sp^{\mu}q^{\nu} && 1\leq i<j\leq n. \notag
\end{align}
Now we replace \eqref{stcond} by
\begin{multline*}
\max\{\abs{t_1},\dots,\abs{t_n},\abs{s_{n+1}},\dots,\abs{s_{n+3}}, \\
\abs{pqS^{-1}t_1^{-1}},\dots,\abs{pqS^{-1}t_n^{-1}}\}<1
<\min\{\abs{s_1},\dots,\abs{s_n}\}.
\end{multline*}
Accordingly, we deform $\T^n$ to $C^n$
such that the set of poles in the interior of $C^n$ is again given
by \eqref{inset}. Then, obviously,
\begin{equation}\label{rhointC}
\kappa^{\A}\int_{C^n} \rho(z;s,t) \d{z}=1.
\end{equation}
Next the integral over $C^n$ is expanded as a sum over
integrals over $\T^{n-m}$ (for the details of such
an expansion, see Section~\ref{secproof}).
Let $N_i$ be a fixed non-negative integer and
assume that $1<\abs{s_iq^{N_i}}<\abs{q}^{-1}$ for $i\in\set{n}$
and $\abs{p}<\min\{\abs{s_1}^{-1},\dots,\abs{s_n}^{-1}\}$.
Then the only poles of the integrand crossing the contour in its deformation
from $C^n$ to $\T^n$ are the poles at
$z_j=s_iq^{\la_i}$ for $\la_i\in\Zet{N_i+1}$, $i\in\set{n}$
and $j\in\set{n+1}$, and the expansion takes the form
\begin{equation}\label{CT}
\int_{C^n} \rho(z;s,t) \d{z}
=\sum_{m=0}^n \sum_{\la^{(m)}}
\int_{\T^{n-m}}\rho_{\la^{(m)}}(z^{(n-m)};s,t) \d{z^{(n-m)}},
\end{equation}
with $z^{(i)}=(z_1,\dots,z_i)$, $\la^{(i)}=(\la_1,\dots,\la_i)\in\Z^i$,
and where the sum over $\la^{(m)}$ is subject to the restriction
that $0\leq \la_i\leq N_i$ for all $i\in\set{m}$. The function
$\rho_{\la^{(m)}}(z^{(n-m)};s,t)$ is obtained from $\rho_{\la^{(0)}}(z^{(n)};s,t)
=\rho(z;s,t)$ by computing the relevant residues.
For the present purposes we only need the explicit form of the kernel
for $m=n$. Writing $\la$ for $\la^{(n)}$ and $\rho_{\la}(z^{(0)};s,t)$ as
$\rho_{\la}(s,t)$, it is given by
\begin{align*}
\rho_{\la}(s,t)&=(2\pi\im)^n (n+1)!\,
\Res_{z_1=s_1q^{\la_1}}\Bigl( \cdots
\Bigl(\Res_{z_n=s_nq^{\la_n}}\Bigl(\frac{\rho(z;s,t)}{z}
\Bigr)\Bigr)\cdots\Bigr) \\
&=\frac{1}{\kappa^{\A}}
\prod_{i=1}^n\biggl(
\frac{\Gamma(S'^{-1}t_i,SS's_i^{-1})}{\Gamma(S'^{-1}s_i^{-1},SS't_i)}
\prod_{j=n+1}^{n+3}
\frac{\Gamma(St_is_j^{-1},s_i^{-1}s_j)}{\Gamma(Ss_i^{-1}s_j^{-1},t_is_j)}
\biggr) \\
&\qquad\times
\prod_{i=1}^n \frac{\theta(S's_iq^{\la_i+\abs{\la}})}{\theta(S's_i)}
\prod_{1\leq i<j\leq n}
\frac{\theta(q^{\la_i-\la_j}s_is_j^{-1})}{\theta(s_is_j^{-1})}
\frac{1}{(qS^{-1}s_is_j)_{\la_i+\la_j}} \\
&\qquad\times\prod_{j=1}^{n+3} \frac{(S's_j)_{\abs{\la}}}
{\prod_{i=1}^n(qs_is_j^{-1})_{\la_i}}
\prod_{i,j=1}^n (s_it_j,qS^{-1}s_it_j^{-1})_{\la_i} \\
&\qquad \times
\prod_{i=1}^n \frac{(SS's_i^{-1})_{\abs{\la}-\la_i} \, q^{i\la_i}}
{(SS't_i,qS't_i^{-1})_{\abs{\la}}}
\end{align*}
for $S'=s_1\cdots s_n$ and $\abs{\la}=\la_1+\cdots+\la_n$.
The other kernels arise as $m$-fold residues and their explicit
form is quite involved. In the remainder we will only use the fact that
$\rho_{\la^{(m)}}(z^{(n-m)};s,t)$ contains the factor
$\prod_{i=1}^{n-m}1/\Gamma(t_is_i)$.

The next step in the computation is to let $t_i$ tend to $q^{-N_i}s_i^{-1}$
in \eqref{CT} for all $i\in\set{n}$. Since for $m\neq n$ the factor
$\prod_{i=1}^{n-m}1/\Gamma(t_is_i)$
vanishes in this limit, the only contribution to the sum over $m$
comes from the term $m=n$. For later reference we state this explicitly;
\begin{equation}\label{limint}
\kappa^{\A}\lim_{\substack{t_i\to q^{-N_i}s_i^{-1}\\ \forall\,i\in\set{n}}}
\int_{C^n} \rho(z;s,t) \d{z}
=\kappa^{\A}\lim_{\substack{t_i\to q^{-N_i}s_i^{-1}\\ \forall\,i\in\set{n}}}
\sum_{\la} \rho_{\la}(s,t),
\end{equation}
with $\la_i$ ranging from $0$ to $N_i$ in the sum on the right.
Thanks to \eqref{rhointC} the left-hand side is equal to $1$,
leading to the following elliptic hypergeometric series identity.
\begin{theorem}\label{thmRos}
For $S=s_1\cdots s_{n+3}$, $S'=s_1\cdots s_n$ we have
\begin{align*}
\sum_{\la}\prod_{i=1}^n&\frac{\theta(S's_iq^{\la_i+\abs{\la}})}{\theta(S's_i)}
\prod_{1\leq i<j\leq n}
\frac{\theta(q^{\la_i-\la_j}s_is_j^{-1})}{\theta(s_is_j^{-1})}
\frac{1}{(qS^{-1}s_is_j)_{\la_i+\la_j}} \\
&\times\prod_{j=1}^{n+3} \frac{(S's_j)_{\abs{\la}}}
{\prod_{i=1}^n(qs_is_j^{-1})_{\la_i}}
\prod_{i,j=1}^n (q^{-N_j}s_is_j^{-1},q^{N_j+1}S^{-1}s_is_j)_{\la_i} \\
&\times \prod_{i=1}^n \frac{(SS's_i^{-1})_{\abs{\la}-\la_i} \, q^{i\la_i}}
{(q^{-N_i}SS's_i^{-1},q^{N_i+1}S's_i)_{\abs{\la}}} \\
&\qquad \qquad \qquad =\prod_{i=1}^n\biggl(
\frac{(qS's_i)_{N_i}}{(qS^{-1}S'^{-1}s_i)_{N_i}}
\prod_{j=n+1}^{n+3} \frac{(qS^{-1}s_is_j)_{N_i}}{(qs_is_j^{-1})_{N_i}}
\biggr),
\end{align*}
where the sum is over $\la=(\la_1,\dots,\la_n)$ such that
$0\leq\la_i\leq N_i$ for all $i\in\set{n}$.
\end{theorem}
The above result is equivalent to the sum proven by Rosengren
in \cite[Corollary 6.3]{Rosengren02}, which is an elliptic version of the
Schlosser's D$_n$ Jackson sum \cite[Theorem 5.6]{Schlosser97}.
In view of our derivation it appears more appropriate to
associate Theorem~\ref{thmRos} with the root system A$_n$.

An alternative way to modify \eqref{stcond} is to take
\begin{equation}\label{stcondp}
\max\{\abs{s_1},\dots,\abs{s_{n+3}},
\abs{pqS^{-1}t_1^{-1}},\dots,\abs{pqS^{-1}t_n^{-1}}\}<1
<\min\{\abs{t_1},\dots,\abs{t_n}\}.
\end{equation}
Again deforming $\T^n$ to $C^n$ so as to let \eqref{inset} be the
set of poles in the interior of $C^n$ we once more
get an integral identity of the form \eqref{rhointC}.
Assuming that $\abs{p}<\min\{\abs{t_1}^{-1},\dots,\abs{t_n}^{-1}\}$
and $1<\abs{t_iq^{N_i}}<\abs{q}^{-1}$ for $i\in\set{n}$, the poles
crossing the contour in its deformation from $C^n$ back to $\T^n$
now correspond to the poles at
$z_j=t_i^{-1}q^{-\la_i}$ for $\la_i\in\{0,\dots,N_i\}$, $i\in\set{n}$
and $j\in\set{n+1}$. Appropriately redefining the integration kernels
the expansion \eqref{CT} still takes place. In particular, this time
\begin{align*}
\rho_{\la}(s,t)&=
(2\pi\im)^n (-1)^n (n+1)! \\
& \qquad \qquad \times
\Res_{z_1=t_1^{-1}q^{-\la_1}}\Bigl( \cdots
\Bigl(\Res_{z_n=t_n^{-1}q^{-\la_n}}\Bigl(\frac{\rho(z;s,t)}{z}
\Bigr)\Bigr)\cdots\Bigr) \\
&=\frac{1}{\kappa^{\A}}
\prod_{i=1}^n \frac{1}{\Gamma(T^{-1}t_i^{-1})}
\prod_{j=1}^{n+3}\biggl(\Gamma(T^{-1}s_j)
\prod_{i=1}^n\Gamma(St_is_j^{-1})\biggr) \\
&\qquad \times\prod_{1\leq i\leq j\leq n}\frac{1}{\Gamma(St_it_j)}
\prod_{1\leq i<j\leq n+3}\frac{1}{\Gamma(Ss_i^{-1}s_j^{-1})} \\
&\qquad \times
\prod_{i=1}^n \frac{\theta(Tt_iq^{\la_i+\abs{\la}})}{\theta(Tt_i)}
\prod_{1\leq i<j\leq n}
\frac{\theta(q^{\la_i-\la_j}t_it_j^{-1})}{\theta(t_it_j^{-1})}
(St_it_j)_{\la_i+\la_j} \\
&\qquad \times
\prod_{j=1}^{n+3} \frac{\prod_{i=1}^n(t_is_j)_{\la_i}}
{(qTs_j^{-1})_{\abs{\la}}}
\prod_{i,j=1}^n
\frac{1}{(qt_it_j^{-1},St_it_j)_{\la_i}} \\
&\qquad \times
\prod_{i=1}^n \frac{(Tt_i,qTS^{-1}t_i^{-1})_{\abs{\la}}\, q^{i\la_i}}
{(qTS^{-1}t_i^{-1})_{\abs{\la}-\la_i}}
\end{align*}
for $T=t_1\cdots t_n$.
Letting $s_i$ tend to $q^{-N_i}t_i^{-1}$ for $i\in\set{n}$ we once
again find that all but the $m=n$ term
vanishes in the sum over $m$ in \eqref{CT}. After the
identification of $(s_{n+1},s_{n+2},s_{n+3})$ with $(b_1,b_2,b_3)$
this yields the following companion to Theorem~\ref{thmRos}.
\begin{theorem}\label{thmeB}
For $T=t_1\cdots t_n$ and $A=qTb_1^{-1}b_2^{-1}b_3^{-1}$ we have
\begin{align*}
\sum_{\la} &
\prod_{i=1}^n \frac{\theta(Tt_iq^{\la_i+\abs{\la}})}{\theta(Tt_i)}
\prod_{1\leq i<j\leq n}
\frac{\theta(q^{\la_i-\la_j}t_it_j^{-1})}{\theta(t_it_j^{-1})}
(q^{1-\abs{N}}A^{-1}t_it_j)_{\la_i+\la_j} \\
&\qquad \times \prod_{i,j=1}^n \frac{(q^{-N_j}t_it_j^{-1})_{\la_i}}
{(qt_it_j^{-1},q^{1-\abs{N}}A^{-1}t_it_j)_{\la_i}}
\prod_{j=1}^3 \frac{\prod_{i=1}^n(t_ib_j)_{\la_i}}
{(qTb_j^{-1})_{\abs{\la}}} \\
&\qquad \times
\prod_{i=1}^n \frac{(Tt_i,q^{\abs{N}}ATt_i^{-1})_{\abs{\la}}
\, q^{i\la_i}}
{(q^{N_i+1}Tt_i)_{\abs{\la}}
(q^{\abs{N}}ATt_i^{-1})_{\abs{\la}-\la_i}} \\
&=\prod_{i,j=1}^n \frac{(At_i^{-1}t_j^{-1})_{\abs{N}-N_i}}
{(At_i^{-1}t_j^{-1})_{\abs{N}}}
\prod_{1\leq i<j\leq n}\frac{(At_i^{-1}t_j^{-1})_{\abs{N}}}
{(At_i^{-1}t_j^{-1})_{\abs{N}-N_i-N_j}} \\
& \qquad \times
\prod_{i=1}^n\prod_{j=1}^3
\frac{(At_i^{-1}b_j)_{\abs{N}}}{(At_i^{-1}b_j)_{\abs{N}-N_i}}
\prod_{i=1}^n(qTt_i)_{N_i}
\prod_{j=1}^3 \frac{1}{(qTb_j^{-1})_{\abs{N}}},
\end{align*}
where the sum is over all $\la=(\la_1,\dots,\la_n)$ such that
$0\leq \la_i\leq N_i$ for each $i\in\set{n}$, and $N=N_1+\cdots+N_n$.
\end{theorem}
The above result corresponds to the elliptic analogue of Bhatnagar's
D$_n$ summation \cite{Bhatnagar99}, and can be transformed into the
identity of Theorem~\ref{thmRos} by changing the summation indices from
$\la_i$ to $N_i-\la_i$ for all $i\in\set{n}$. Actually, the described
residue calculus with the simplest choices $N_i=0$ (i.e., when there
remains only one trivial term in the sum of Theorem~\ref{thmeB}) will
be used in Section~\ref{secaltproof} in the alternative proof of
Theorem~\ref{thmbetaA}. The full sums of Theorems~\ref{thmRos}
and \ref{thmeB} then follow from the application of general
residue calculus.

\subsection{Proof of Theorem~\ref{inversionthmAC}}\label{secproof}

We begin by introducing some notation and definitions.
For $a\in\Zet{n+1}$ let $Z_0=1$, $Z_a=z_1\cdots z_a$,
$\overline{W}_0=1$, $\overline{W}_a=w_{n-a+1}\cdots w_n$,
$z^{(a)}=(z_1,\dots,z_a)$ and $w^{(a)}=(w_1,\dots,w_a)$.
Note that $Z_n=z_{n+1}^{-1}$, $\overline{W}_n=w_{n+1}^{-1}$,
$z^{(n)}=z$ and $w^{(n)}=w$.
Dropping the superscript $(\A,\C)$ in $\Delta(z,w,x;t)$ we recursively define
\begin{subequations}\label{res1}
\begin{align}\label{res1a}
\Delta(z^{(n-a)},w,x;t)
&=\Res_{z_{n-a+1}=t^{-1}w_{n-a+1}}
\frac{\Delta(z^{(n-a+1)},w,x;t)}{z_{n-a+1}} \\
&=-\Res_{z_{n-a+1}=t^a \overline{W}_a^{-1}Z_{n-a}^{-1}}
\frac{\Delta(z^{(n-a+1)},w,x;t)}{z_{n-a+1}}.
\label{res1b}
\end{align}
\end{subequations}
for $a\in\set{n}$.
The equality of the two expressions on the right easily follows from the
A$_n$ symmetry of $\Delta(z,w,x;t)$ in the $z$-variables.
Indeed, the above recursions imply that
$\Delta(z^{(n-a)},w,x;t)=\Delta(\sigma(z^{(n-a)}),w,x;t)$ for $\sigma\in S_{n-a}$
and that $\Delta(z^{(n-a)},w,x;t)$ is invariant under the variable change
$z_i\to t^a \overline{W}_a^{-1}Z_{n-a}^{-1}$
(and hence $Z_{n-a}\to t^a \overline{W}_a^{-1}z_i^{-1}$)
for arbitrary fixed $i\in\set{n-a}$.
These two symmetries of course generate a group of dimension $(n-a+1)!$
isomorphic to $S_{n-a+1}$, and for $a=0$ correspond to the A$_n$ symmetry of
$\Delta(z,w,x;t)$.
By the C$_n$ symmetry of $\Delta(z,w,x;t)$ in the $w$-variables it also
follows that $\Delta(z^{(n-a)},w,x;t)$ has C$_{n-a}$ symmetry in the
variables $(w_1,\dots,w_{n-a})$ and
C$_a$ symmetry in the variables $(w_{n-a+1},\dots,w_n)$.
Hence for $k\in\set{n-a+1}$ and $\sigma\in\{-1,1\}$
\begin{multline}\label{res2}
\Res_{z_{n-a+1}=t^{-1}w_k^{\sigma}} \frac{\Delta(z^{(n-a+1)},w,x;t)}{z_{n-a+1}} \\
=-\Res_{z_{n-a+1}=t^a w_k^{-\sigma}\overline{W}_{a-1}^{-1}Z_{n-a}^{-1}}
\frac{\Delta(z^{(n-a+1)},w,x;t)}{z_{n-a+1}} \\
=\Delta(z^{(n-a)},w,x;t)|_{w_{n-a+1}\leftrightarrow w_k^{\sigma}}.
\end{multline}

Using \eqref{Glim} and induction, the explicit
form for $\Delta(z^{(n-a)},w,x;t)$ is easily found to be
\begin{align}\label{delta}
\Delta(z^{(n-a)},w,x;t)&=
\frac{1}{\pq{a}}
\frac{\prod_{i=1}^n\prod_{j=1}^{n+1}
\Gamma(tw_i^{\pm}x_j)}
{\prod_{i=1}^n\Gamma(w_i^{\pm 2})
\prod_{1\leq i<j\leq n}\Gamma(w_i^{\pm}w_j^{\pm})} \\
& \quad \times \prod_{i=n-a+1}^n \biggl[\Gamma(w_i^{-2})
\prod_{j=1}^{n-a} \Gamma(w_i^{-1}w_j^{\pm})\biggr]
\prod_{n-a+1\leq i<j\leq n}\! \frac{\Gamma(w_i^{-1}w_j^{-1})}{\Gamma(w_iw_j)}
\notag \\
&\quad \times \prod_{j=1}^{n-a+1}
\frac{\prod_{i=1}^{n-a} \Gamma(t^{-1}w_i^{\pm}z_j^{-1})]}
{\prod_{i=n-a+1}^n \Gamma(tw_i^{\pm}z_j)} \notag \\
&\quad \times\prod_{1\leq i<j\leq n-a+1}
\frac{1}{\Gamma(z_iz_j^{-1},z_i^{-1}z_j,t^2z_iz_j,t^{-2}z_i^{-1}z_j^{-1})},
\notag
\end{align}
where, to keep the expression from spilling over, we have set
$z_{n-a+1}:=t^a\overline{W}_a^{-1} Z_{n-a}^{-1}$.
Note that this also makes all of the claimed symmetries of
$\Delta(z^{(n-a)},w,x;t)$ manifest.

After these preliminaries we can state our first intermediate result.
Let
\begin{equation}\label{defI}
I(x;t):=\int_{\T^n}\int_{C^n_w}\Delta(z,w,x;t)f(z) \d{z}\d{w}.
\end{equation}
\begin{proposition}\label{propexp}
There holds
\begin{multline}\label{eqpropexp}
I(x;t)=\sum_{a=0}^n (4\pi\im)^a \binom{n}{a}\frac{(n+1)!}{(n-a+1)!}
\int_{\T^n}\int_{\T^{n-a}} \Delta(z^{(n-a)},w,x;t)  \\
\times f(z^{(n-a)},t^{-1}w_{n-a+1},\dots,t^{-1}w_n)
\d{z^{(n-a)}}\d{w}.
\end{multline}
\end{proposition}

We break up the proof into two lemmas, the first of which requires
some more notation.
Write $\An$ of \eqref{An} as $\An^n$ and, more generally, define
$\An^{n-a}$ for $a\in\Zet{n+1}$ as
\begin{multline}\label{Ana}
\An^{n-a}=\{z\in\Complex^{n-a}|~\abs{t}^n-\epsilon
<\abs{z_j}<\abs{t}^{-1}+\epsilon,~
j\in\set{n-a}, \\
\text{and}~\abs{Z_{n-a}^{-1}}<\abs{t}^{-a-1}+\epsilon\}
\end{multline}
for infinitesimally small but positive $\epsilon$.
For $a=0$ the condition $\abs{t}^n-\epsilon <\abs{z_j}$
becomes superfluous and we recover \eqref{An}.
We will sometimes somewhat loosely say that $z_j=u$ is not in $\An^{n-a}$
when we really mean that $z^{(n-a)}=(z_1,\dots,z_{n-a})$ with $z_j=u$
is not in $\An^{n-a}$.
Since, $\abs{Z_{n-a-1}^{-1}}=\abs{Z_{n-a}^{-1}}\abs{z_{n-a}}
\leq \abs{t}^{-a-2}+\epsilon$ we have the
filtration
$\An^0\subset \An^1\subset \dots \subset \An^{n-a}$,
so that $f(z^{(n-a)},t^{-1}w_{n-a+1},\dots,t^{-1}w_n)$
is holomorphic on $\An^{n-a}$.

We must also generalize the definition of $C^n_w$ given in
Theorem~\ref{inversionthmAC}. To simplify notations
we drop the explicit $w$ dependence of $C^n_w$ and
write $C^n_w$ as $C^n_0$, with $0$ a new label (not
related to the $w$-variables). More generally we
define $C^{n-a}_m\subset\An^{n-a}$ for
$m\in\Zet{n}$ and $a\in\Zet{m+1}$ as deformations of
$\T^{n-a}$ such that for
fixed $w\in\Complex^{n-a}$ and $i\in\set{n-a}$,
\begin{subequations}\label{inex1}
\begin{align}
z_j&=t^{-1}w_i^{\pm} \;
\text{lies in the interior of $C^{n-a}_m$ for $j\in\{m-a+1,\dots,n-a\}$} \\
z_j&=t^{-1}w_i^{\pm}\;
\text{lies in the exterior of $C^{n-a}_m$ for $j\in\set{m-a}$} \\
Z_{n-a}^{-1}&=t^{-a-1}w_i^{\pm} \overline{W}_a\;
\text{lies in the exterior of $C^{n-a}_m$.}
\end{align}
\end{subequations}
For $m=a=0$ this definition simplifies to \eqref{poles2}.
If $C^{n-a}_m$ and $\hat{C}^{n-a}_m$  both satisfy \eqref{inex1}
they will be referred to as homotopic.

Finally, we introduce the shorthand notation
\begin{equation*}
\F(z^{(n-a)},w):=
\Delta(z^{(n-a)},w,x;t) f(z^{(n-a)},t^{-1}w_{n-a+1},\dots,t^{-1}w_n),
\end{equation*}
where the $x$ and $t$ dependence of $\F$ have been suppressed.

\begin{lemma}\label{lemexp}
For $m\in\Zet{n}$
\begin{equation}\label{eqlemexp}
I(x;t)=\sum_{a=0}^m (4\pi\im)^a \binom{m}{a}\frac{n!}{(n-a)!}
\int_{\T^n}\int_{C^{n-a}_m}
\F(z^{(n-a)},w)
\d{z^{(n-a)}}\d{w}.
\end{equation}
\end{lemma}

\begin{proof}
We prove this by induction on $m$.
Since for $m=0$ we recover definition \eqref{defI} of $I(x;t)$,
we only need to establish the induction step.
Writing the expression on the right of \eqref{eqlemexp} for fixed $m$ as
$I_m(x;t)$, the problem is to show that
$I_m(x;t)=I_{m+1}(x;t)$ for $m\in\Zet{n-1}$.

Since
\begin{equation}\label{van}
\Delta(z^{(n-a)},w,x;t)=0\;\text{if $w_i=w_j$ for $1\leq i<j\leq n-a$},
\end{equation}
we may without loss of generality assume that all components of $w^{(n-a)}$
are distinct when considering the $z^{(n-a)}$-integration for fixed $w$.

In the integral over $z^{(n-a)}$ we deform $C^{n-a}_m$ to
$B^{n-a}_m\subset\An^{n-a}$ such that
\begin{align*}
z_j&=t^{-1}w_i^{\pm} \;
\text{lies in the interior of $B^{n-a}_m$ for $j\in\{m-a+1,\dots,n-a-1\}$} \\
z_j&=t^{-1}w_i^{\pm}\;
\text{lies in the exterior of $B^{n-a}_m$ for $j\in\set{m-a}\cup\{n-a\}$} \\
Z_{n-a}^{-1}&=t^{-a-1}w_i^{\pm} \overline{W}_a\;
\text{lies in the exterior of $B^{n-a}_m$.}
\end{align*}
To see how this deformation changes the integral over $z^{(n-a)}$ we
need to investigate the location of the poles of the integrand.
{}From \eqref{ellGamma} and \eqref{delta} it follows that
$\Delta(z^{(n-a)},w,x;t)$ has poles at
\begin{align*}
z_j&=t^{-1}w_i^{\pm} p^{\mu}q^{\nu}&&\text{for}~i,j\in\set{n-a}\\
Z_{n-a}^{-1}&=t^{-a-1}w_i^{\pm}\overline{W}_a p^{\mu}q^{\nu}&&
\text{for}~i\in\set{n-a}\\
z_j&=t^{-1}w_i^{\pm} p^{\mu+1}q^{\nu+1}&&\text{for}~
i\in\{n-a+1,\dots,n\},~j\in\set{n-a}\\
Z_{n-a}^{-1}&=t^{-a-1}w_i^{\pm}\overline{W}_a p^{\mu+1}q^{\nu+1}
\hspace{-5mm}&&
\text{for}~i\in\{n-a+1,\dots,n\}.
\end{align*}
We note in particular that the terms in the last line of \eqref{delta}
do not imply any poles for $\Delta(z^{(n-a)},w,x;t)$ thanks to the
reflection equation \eqref{reflex}.
Of the poles listed above only those corresponding to the first
two lines with $(\mu,\nu)=(0,0)$ are in $\An^{n-a}$,
i.e., the poles at $z_j=t^{-1}w_i^{\pm}$ for $i,j\in\set{n-a}$
and at $Z_{n-a}^{-1}=t^{-a-1}w_i^{\pm}\overline{W}_a$ for
$i\in\set{n-a}$.
Indeed, for $(\mu,\nu)\neq (0,0)$, all of the above $z_j$ satisfy
(since $w\in\T^n$) $\abs{z_j}\leq \abs{t}^{-1}M<\abs{t}^n$,
incompatible with \eqref{Ana}.
Likewise, for $(\mu,\nu)\neq (0,0)$, the above listed poles
for $Z_{n-a}^{-1}$ satisfy $\abs{Z_{n-a}^{-1}}\leq \abs{t}^{-a-1}M<\abs{t}^{n-a}$.
But from \eqref{Ana} we have
$\abs{z_j}<\abs{t}^{-1}$ for all $j\in\set{n-a}$ which implies that
$\abs{Z_{n-a}^{-1}}>\abs{t}^{n-a}$.

Comparing the definitions $C^{n-a}_m$ and $B^{n-a}_m$ we thus see that
the only difference between the integral over the former and the latter
is that the poles at $z_{n-a}=t^{-1}w_k^{\pm}$ for $k\in\set{n-a}$
have moved to the exterior of $B^{n-a}_m$
To compensate for this discrepancy we need to calculate the
residues (denoted by $R_{k,\pm}$) of the integrand at
$z_{n-a}=t^{-1}w_k^{\pm}$ and integrate this over
$B^{n-a-1}_{m;k,\pm}\subset\An^{n-a-1}$.
Here $B^{n-a-1}_{m;k,\pm}$ corresponds to ``$C^{n-a}_m$ restricted to
$z_{n-a}=t^{-1}w_k^{\pm}$.'' That is,
for fixed $w\in\T^{n-a}$ and $i\in\set{n-a-1}$, $i\neq k$,
\begin{align*}
z_j&=t^{-1}w_i^{\pm} \;
\text{lies in the interior of $B^{n-a-1}_{m;k,\sigma}$ }
\\ & \qquad \qquad \qquad \qquad \qquad
\text{for $j\in\{m-a+1,\dots,n-a-1\}$}
\\ z_j&=t^{-1}w_i^{\pm} \;
\text{lies in the exterior of $B^{n-a-1}_{m;k,\sigma}$
for $j\in\set{m-a}$} \\
Z_{n-a-1}^{-1}&=t^{-a-2}w_i^{\pm} w_k^{\sigma}\overline{W}_ap^{\mu} q^{\nu}\;
\text{lies in the exterior of $B^{n-a-1}_{m;k,\sigma}$.}
\end{align*}
The exclusion of $i=k$ is justified by the fact that
$R_{k,\sigma}$ is free of poles at the above when $i=k$
thanks to \eqref{van}.
By \eqref{res2} with $a\to a+1$ and the fact that $f$
is holomorphic on $\An^{n-a}$ we get
\begin{multline*}
I_m(x;t)=\sum_{a=0}^m (4\pi\im)^a \binom{m}{a}\frac{n!}{(n-a)!}
\biggl[\int_{\T^n}\int_{B^{n-a}_m} \F(z^{(n-a)},w)
\d{z^{(n-a)}}\d{w}\\
+2\pi\im \sum_{k=1}^{n-a}\sum_{\sigma\in\{\pm 1\}}
\int_{\T^n}\int_{B^{n-a-1}_{m;k,\sigma}}
\F(z^{(n-a-1)},w)|_{w_{n-a}\leftrightarrow w_k^{\sigma}}
\d{z^{(n-a-1)}}\d{w}\biggr].
\end{multline*}
We now change integration variables in both terms inside the square brackets.
In the first term we substitute $z_{n-a}\leftrightarrow z_{m-a+1}$ and in the
second term (or rather its summand for fixed $\sigma$)
we substitute $w_k^{\sigma}\leftrightarrow w_{n-a}$.
Using the permutation symmetry of $f$ and noting that
$B^{n-a}_m|_{z_{n-a}\leftrightarrow z_{m-a+1}}$ is homotopic to
$C^{n-a}_{m+1}$ and $B^{n-a-1}_{m;k,\sigma}|_{w_k\leftrightarrow w_{n-a}}=
B^{n-a-1}_{m;n-a,1}$ is homotopic to $C^{n-a-1}_{m+1}$,
this yields
\begin{multline*}
I_m(x;t)=\sum_{a=0}^m (4\pi\im)^a \binom{m}{a}\frac{n!}{(n-a)!}
\biggl[\int_{\T^n} \int_{C^{n-a}_{m+1}} \F(z^{(n-a)},w)
\d{z^{(n-a)}}\d{w}\\
+(4\pi\im)(n-a)
\int_{\T^n} \int_{C^{n-a-1}_{m+1}} \F(z^{(n-a-1)},w)
\d{z^{(n-a-1)}}\d{w}\biggr].
\end{multline*}
Finally changing the summation index $a\to a-1$ in the sum
corresponding to the second term inside the square brackets,
and using the standard binomial recursion leads to the desired
$I_m(x;t)=I_{m+1}(x;t)$.
\end{proof}

{}From Lemma~\ref{lemexp} with $m=n-1$ it follows that
\begin{equation*}
I(x;t)=\sum_{a=0}^{n-1} (4\pi\im)^a \binom{n-1}{a}\frac{n!}{(n-a)!}
\int_{\T^n}\int_{C^{n-a}} \F(z^{(n-a)},w)
\d{z^{(n-a)}}\d{w},
\end{equation*}
where the label $n-1$ of $C^{n-a}_{n-1}$ has been dropped,
having served its purpose.

The second lemma needed to prove Proposition~\ref{propexp} should
thus read as follows.
\begin{lemma}
There holds
\begin{multline}\label{lemincex}
\sum_{a=0}^{n-1} (4\pi\im)^a \binom{n-1}{a}\frac{n!}{(n-a)!}
\int_{\T^n}\int_{C^{n-a}} \F(z^{(n-a)},w)
\d{z^{(n-a)}}\d{w} \\
=\sum_{a=0}^n (4\pi\im)^a \binom{n}{a}\frac{(n+1)!}{(n-a+1)!}
\int_{\T^n}\int_{\T^{n-a}} \F(z^{(n-a)},w)
\d{z^{(n-a)}}\d{w}.
\end{multline}
\end{lemma}

\begin{proof}
The only difference between the two integrals over $z^{(n-a)}$
in \eqref{lemincex} is that $C^{n-a}$ --- given by \eqref{inex1}
with $m=n-1$ --- has the poles of the the integrand at
$z_{n-a}=t^{-1}w_k^{\pm}$ for $k\in\set{n-a}$ in its interior
and the poles at $Z_{n-a}^{-1}=t^{-a-1}w_k^{\pm} \overline{W}_a$
for $k\in\set{n-a}$
(i.e., $z_{n-a}=t^{a+1}w_k^{\pm} \overline{W}_a^{-1}Z_{n-a-1}^{-1}$)
in its exterior, whereas $\T^{n-a}$ has the latter in its interior
and the former in its exterior.
Hence, applying \eqref{res2} and
\begin{multline*}
f(z^{(n-a-1)},t^{a+1}w_k^{-\sigma}\overline{W}_a^{-1}Z_{n-a-1}^{-1},
t^{-1}w_{n-a+1},\dots,t^{-1}w_n) \\
=f(z^{(n-a-1)},t^{-1}w_k^{\sigma},t^{-1}w_{n-a+1},\dots,t^{-1}w_n)
\end{multline*}
as follows from the A$_n$ symmetry of $f$, we get
\begin{multline*}
\int_{\T^n}\int_{C^{n-a}} \F(z^{(n-a)},w)
\d{z^{(n-a)}}\d{w}
\to
\int_{\T^n}\int_{\T^{n-a}} \F(z^{(n-a)},w)
\d{z^{(n-a)}}\d{w} \\
+4\pi\im
\sum_{k=1}^{n-a}\sum_{\sigma\in\{\pm 1\}}
\int_{\T^n}\int_{\T^{n-a-1}}
\F(z^{(n-a-1)},w)|_{w_{n-a}\leftrightarrow w_k^{\sigma}}
\d{z^{(n-a-1)}}\d{w} \\
=\int_{\T^n}\int_{\T^{n-a}} \F(z^{(n-a)},w)
\d{z^{(n-a)}}\d{w} \\
+8\pi\im(n-a)
\int_{\T^n}\int_{\T^{n-a-1}} \F(z^{(n-a-1)},w)
\d{z^{(n-a-1)}}\d{w}.
\end{multline*}
Here the last expression on the right follows after the
variable change $w_k^{\sigma}\leftrightarrow w_{n-a}$ in the
second double integral.
We wish to emphasize that one of the factors $2$ in $4\pi\im\sum_k \dots $ is due
to the fact that the poles at $z_{n-a}=t^{-1}w_k^{\pm}$ and
$z_{n-a}=t^{a+1}w_k^{\mp} \overline{W}_a^{-1}Z_{n-a-1}^{-1}$
yields the same contribution by virtue of \eqref{res2}.
In what follows we further examine the
contributions arising from $z_{n-a}=t^{-1}w_i^{\pm}$.

The reason for putting an arrow instead of an equal sign in the above
is that the expression on the right is overcounting poles and needs an
additional correction term.
Indeed, we have computed the residues of $\F(z^{(n-a)},w)/z_{n-a}$ at
its poles $z_{n-a}=t^{-1}w_k^{\pm}$. By exploiting the symmetry \eqref{res2}
and by making a variable change in the $w$-variables this effectively boiled
down to picking up the residue at $z_{n-a}=t^{-1}w_{n-a}$ exactly $2(n-a)$ times.
This residue, given by $\F(z^{(n-a-1)},w)$, has poles at
\begin{equation}\label{wrongpoles}
Z_{n-a-1}=t^{a+2}w_i^{-\sigma}\overline{W}_{a+1}^{-1},\quad k\in\set{n-a-1}.
\end{equation}
These poles are in the interior of
$\T^{n-a-1}$ and thus contribute to the integral over
$\F(z^{(n-a-1)},w)$.
But \eqref{wrongpoles} times $z_{n-a}=t^{-1}w_{n-a}$ yields
$Z_{n-a}=t^{a+1}w_i^{-\sigma}\overline{W}_a^{-1}$, or, equivalently,
$Z_{n-a}^{-1}=t^{-a-1}w_i^{\sigma}\overline{W}_a$.
According to \eqref{inex1} with $m=n-1$ these poles lie in the exterior
of $C^{n-a}$ and
hence the poles \eqref{wrongpoles} should not be contributing at all!
Consequently we need to subtract the further term
\begin{multline*}
-2(2\pi\im)^2(n-a)\sum_{k=1}^{n-a-1}
\sum_{\sigma\in\{\pm 1\}}  \\ \qquad \qquad \times
\int_{\T^n}\int_{\T^{n-a-2}}
\F(z^{(n-a-2)},w)|_{w_{n-a-1}\leftrightarrow w_k^{\sigma}}
\d{z^{(n-a-2)}}\d{w} \\
=-(4\pi\im)^2(n-a)(n-a-1)\int_{\T^n}\int_{\T^{n-a-2}}
\F(z^{(n-a-2)},w)
\d{z^{(n-a-2)}}\d{w},
\end{multline*}
where we have used the second equality in \eqref{res2} with $a\to a+2$, and
the A$_n$ symmetry of $f$. Therefore
\begin{multline*}
\int_{\T^n}\int_{C^{n-a}} \F(z^{(n-a)},w)
\d{z^{(n-a)}}\d{w}
=\int_{\T^n}\int_{\T^{n-a}} \F(z^{(n-a)},w)
\d{z^{(n-a)}}\d{w} \\
+2(4\pi\im)(n-a)
\int_{\T^n}\int_{\T^{n-a-1}} \F(z^{(n-a-1)},w)
\d{z^{(n-a-1)}}\d{w} \\
+(4\pi\im)^2(n-a)(n-a-1)\int_{\T^n}\int_{\T^{n-a-2}}
\F(z^{(n-a-2)},w)
\d{z^{(n-a-2)}}\d{w}.
\end{multline*}
Substituting this in the left-hand side of \eqref{lemincex}, shifting
$a\to a-1$ and $a\to a-2$ in the sums corresponding to the integrals
over $\T^{n-a-1}$ and $\T^{n-a-2}$ and using the binomial identity
\begin{equation*}
\binom{n-1}{a}+2\binom{n-1}{a-1}+\binom{n-1}{a-2}
=\binom{n+1}{a}=\frac{n+1}{n-a+1}\binom{n}{a}
\end{equation*}
yields the wanted right-hand side of \eqref{lemincex}, completing the proof.
\end{proof}

In the integral on the right-hand side of \eqref{eqpropexp}
we make the variable changes $z_i\to z_{i+a}$ for $i\in\set{n-a}$
followed by $t^{-1}w_{n-a+i}\to z_i$ for $i\in\set{a}$.
By the permutation symmetry of $f$ this gives
\begin{multline}\label{cont}
I(x;t)=\sum_{a=0}^n (4\pi\im)^a \binom{n}{a}\frac{(n+1)!}{(n-a+1)!} \\
\times
\int_{\T^{n-a}} \int_{(t^{-1}\T)^a\times \T^{n-a}}
\Delta(z,w^{(n-a)},x;t)f(z) \d{z}\d{w^{(n-a)}},
\end{multline}
with
\begin{equation*}
\Delta(z,w^{(n-a)},x;t):=
\Delta((z_{a+1},\dots,z_n),(w_1,\dots,w_{n-a},tz_1,\dots,tz_a),x;t)
\end{equation*}
given by the somewhat unwieldy expression
\begin{align}\label{unwieldy}
&\Delta(z,w^{(n-a)},x;t)=
\frac{1}{\pq{a}}
\prod_{1\leq i<j\leq n-a}\frac{1}{\Gamma(w_i^{\pm}w_j^{\pm})} \\
&\qquad\times
\prod_{j=1}^a \frac{\prod_{i=1}^{n+1}\Gamma(t^2x_iz_j,x_iz_j^{-1})}
{\Gamma(t^2z_j^2)
\prod_{i=a+1}^{n+1}\Gamma(z_iz_j^{-1},t^2z_iz_j)} \notag \\
&\qquad\times\prod_{i=1}^{n-a} \frac{
\prod_{j=1}^{n+1}\Gamma(tw_i^{\pm}x_j)
\prod_{j=a+1}^{n+1} \Gamma(t^{-1}w_i^{\pm}z_j^{-1})}
{\Gamma(w_i^{\pm 2})\prod_{j=1}^a\Gamma(tw_i^{\pm}z_j)} \notag \\
&\qquad\times
\prod_{1\leq i<j\leq a}
\frac{1}{\Gamma(z_iz_j^{-1},z_i^{-1}z_j,t^2z_iz_j,t^2z_iz_j)} \notag \\
&\qquad\times
\prod_{a+1\leq i<j\leq n+1} \frac{1}
{\Gamma(z_iz_j^{-1},z_i^{-1}z_j,t^2z_iz_j,t^{-2}z_i^{-1}z_j^{-1})}. \notag
\end{align}

It is easily checked that the only poles of
$\Delta(z,w^{(n-a)},x;t)$ in the $z_j$-plane located on the annulus
$1<\abs{z_j}<\abs{t}^{-1}$ for
$j\in\set{a}$ are given by $z_j=x_i$ for $i\in\set{n}$.
For example, the pole at $z_{n+1}=t^{-1}w_i^{\pm} p^{\mu}q^{\nu}$
corresponds to a pole in the $z_j$-plane ($j\in\set{a}$)
at $z_j=tw_i^{\pm}(\prod_{k=1;k\neq j}^n z_k)^{-1}p^{-\mu}q^{-\nu}$
with absolute value $\abs{z_j}=\abs{t}^a\abs{p}^{-\mu}\abs{q}^{-\nu}$.
For $(\mu,\nu)=(0,0)$ this implies $\abs{z_j}<1$ and for
$(\mu,\nu)\neq(0,0)$ this implies $\abs{z_j}>\abs{t}^{a}M^{-1}>\abs{t}^{n-1}M^{-1}>
\abs{t}^{-1}$.
As another example, the pole at $z_j^2=t^{-2}p^{\mu+1}q^{\nu+1}$ has
absolute value
$\abs{z_j}^2\leq \abs{t^{-2}pq}\leq \abs{t^{-2}}M^2\leq \abs{t}^{2n}
\leq 1$, et cetera.
Consequently, in deflating the contours $t^{-1}\T$ to $\T$
the poles at $z_j=x_i$ for $j\in\set{a}$ and $i\in\set{n}$ move
from the interior to the exterior but no other poles of the integrand
cross the contours of integration. Recursively defining the necessary residues as
\begin{multline}\label{recdef}
\Delta^{k_1,\dots,k_b}(z^{(n-b)},w^{(n-a)},x;t) \\
=\Res_{\zeta=x_{k_1}}
\frac{\Delta^{k_2,\dots,k_b}((z_1,\dots,z_{a-b},\zeta,
z_{a-b+2},\dots,z_{n-b}),w^{(n-a)},x;t)}{\zeta} \\
\end{multline}
for $b\in\set{a}$ and $k_1,\dots,k_b\in\set{r}$ with $k_i\neq k_j$ we
get our third lemma.
\begin{lemma}\label{lemexp2}
There holds
\begin{multline*}
I(x;t)=\sum_{a=0}^n\sum_{b=0}^a
\frac{2^a n!(n+1)!(2\pi\im)^{a+b} }{(a-b)!(n-a)!(n-a+1)!}
\sum_{1\leq k_1<\dots<k_b\leq n}
\int_{\T^{n-a}} \int_{\T^{n-b}}
\\ \times
\Delta^{k_1,\dots,k_b}(z^{(n-b)},w^{(n-a)},x;t)
f(z^{(n-b)},x_{k_1},\dots,x_{k_b}) \d{z^{(n-b)}}\d{w^{(n-a)}}.
\end{multline*}
\end{lemma}

\begin{proof}
Since $\Delta(z,w^{(n-a)},x;t)$ exhibits permutation symmetry in the
variables $z_1,$ $\ldots,z_a$ it follows that
\begin{equation*}
\Delta^{k_1,\dots,k_b}(z^{(n-b)},w^{(n-a)},x;t)
=\Delta^{\varw(k_1,\dots,k_b)}(z^{(n-b)},w^{(n-a)},x;t)
\end{equation*}
for $\varw\in S_a$.
Hence $\Delta^{k_1,\dots,k_b}(z^{(n-b)},w^{(n-a)},x;t)$ is invariant under
permutation of the variables $x_{k_1},\dots,x_{k_b}$.

After this preliminary comment we will show
by induction on $c$ that
\begin{multline}\label{expabc}
\int_{(t^{-1}\T)^a\times\T^{n-a}} \Delta(z,w^{(n-a)},x;t)f(z) \d{z} \\
=\sum_{b=0}^c (2\pi\im)^b \binom{c}{b}
\sum_{\substack{k_1,\dots,k_b=1 \\ k_i\neq k_j}}^n
\int_{(t^{-1}\T)^{a-c}\times \T^{n-a+c-b}} \qquad \qquad \\
\times
\Delta^{k_1,\dots,k_b}(z^{(n-b)},w^{(n-a)},x;t)
f(z^{(n-b)},x_{k_1},\dots,x_{k_b}) \d{z^{(n-b)}}
\end{multline}
for $c\in\Zet{a+1}$. Since this is trivially true for $c=0$
we only need to establish the induction step.
Let $L(a,c)$ denote the left hand side of \eqref{expabc}. Then
\begin{align*}
L(a,c)&=\sum_{b=0}^c (2\pi\im)^b \binom{c}{b}
\sum_{\substack{k_1,\dots,k_b=1 \\ k_i\neq k_j}}^n
\int_{(t^{-1}\T)^{a-c-1}\times\T^{n-a+c-b+1}}  \\
& \qquad \times \Delta^{k_1,\dots,k_b}(z^{(n-b)},w^{(n-a)},x;t)
f(z^{(n-b)},x_{k_1},\dots,x_{k_b}) \d{z^{(n-b)}} \\
&\quad +\sum_{b=0}^c (2\pi\im)^{b+1} \binom{c}{b}
\sum_{\substack{k_0,k_1,\dots,k_b=1 \\ k_i\neq k_j}}^n
\int_{(t^{-1}\T)^{a-c-1}\times\T^{n-a+c-b}}  \\
& \qquad \times
\Delta^{k_0,\dots,k_b}(z^{(n-b-1)},w^{(n-a)},x;t)
f(z^{(n-b-1)},x_{k_0},\dots,x_{k_b}) \d{z^{(n-b-1)}}.
\end{align*}
In the seccond sum we shift the summation index $b\to b+1$ and relabel the $k_i$
as $k_{i+1}$. By the standard binomial recurrence we then find that
$L(a,c)=L(a,c+1)$ as desired.

To now obtain the expansion of Lemma~\ref{lemexp2} we choose $c=a$ in \eqref{expabc}
and use the permutation symmetry of the integrand in the $x_{k_i}$ to
justify the simplification
\begin{equation*}
\sum_{\substack{k_1,\dots,k_b=1 \\ k_i\neq k_j}}^n (\dots)\to
b!\sum_{1\leq k_1<\dots<k_b\leq n}(\dots).
\end{equation*}
Substituting the resulting expression in \eqref{cont} completes the proof.
\end{proof}

By changing the order of the sums as well as the order of the integrals,
the identity of Lemma~\ref{lemexp2} can be rewritten as
\begin{multline*}
I(x;t)=\sum_{b=0}^n \sum_{1\leq k_1<\dots<k_b\leq n}
\sum_{a=b}^n
\frac{2^an!(n+1)!(2\pi\im)^{a+b}}{(a-b)!(n-a)!(n-a+1)!} \\
\times
\int_{\T^{n-b}}
\biggl[\int_{\T^{n-a}} \Delta^{k_1,\dots,k_b}(z^{(n-b)},w^{(n-a)},x;t)
\d{w^{(n-a)}}\biggr] \\
\times
f(z^{(n-b)},x_{k_1},\dots,x_{k_b}) \d{z^{(n-b)}}.
\end{multline*}

In what may well be considered the second part of the proof of
Theorem~\ref{inversionthmAC} we will 
prove the following proposition.
\begin{proposition}\label{propnul}
For $b\in\Zet{n}$ and fixed $1\leq k_1<\dots<k_b\leq n$ we have
\begin{multline}\label{nul}
\sum_{a=b}^n
\frac{(4\pi\im)^a}{(n-a+1)!} \binom{n-b}{n-a}
\int_{\T^{n-b}}
\biggl[\int_{\T^{n-a}} \Delta^{k_1,\dots,k_b}(z^{(n-b)},w^{(n-a)},x;t)
\d{w^{(n-a)}}\biggr] \\
\times
f(z^{(n-b)},x_{k_1},\dots,x_{k_b}) \d{z^{(n-b)}}=0.
\end{multline}
\end{proposition}
{}From this result it follows that the only non-vanishing contribution to
$I(x;t)$ comes from $b=a=n$ and $(k_1,\dots,k_n)=(1,\dots,n)$.
As we shall see shortly,
\begin{equation}\label{wantedterm}
\Delta(z^{(0)},w^{(0)},x;t):=
\Delta^{1,\dots,n}(z^{(0)},w^{(0)},x;t)=\frac{1}{\pq{2n}},
\end{equation}
leading to
$\kappa^{\A}\kappa^{\C} I(x;t)=f(x)$
as claimed by the theorem.

\begin{proof}[Proof of Proposition~\ref{propnul}]
It is sufficient to prove the proposition for
$(k_1,\dots,k_b)=(1,\dots,b)$. Other choices of the $k_i$ (for fixed $b$)
simply follow by a relabelling of the $x_i$ variables. The advantage of this
particular choice of $k_i$ is that many of the expressions below
significantly simplify.

The first ingredient needed for the proof is the actual computation of
\begin{equation*}
\Delta(z^{(n-b)},w^{(n-a)},x;t):=\Delta^{1,\dots,b}(z^{(n-b)},w^{(n-a)},x;t).
\end{equation*}
{}From definition \eqref{recdef}, and the equations \eqref{unwieldy} and \eqref{Glim}
it is not hard to show that
\begin{align*}
\Delta&(z^{(n-b)},w^{(n-a)},x;t)=\frac{1}{\pq{a+b}}
\prod_{1\leq i<j\leq n-a}\frac{1}{\Gamma(w_i^{\pm}w_j^{\pm})} \\
&\times
\prod_{i=1}^b
\frac{\prod_{j=b+1}^{n+1}\Gamma(x_i^{-1}x_j,t^2x_ix_j)}
{\prod_{j=1}^{n-b+1} \Gamma(x_i^{-1}z_j,t^2x_iz_j)}
\prod_{j=1}^{a-b} \frac{
\prod_{i=b+1}^{n+1}\Gamma(t^2x_iz_j,x_iz_j^{-1})}{\Gamma(t^2z_j^2)
\prod_{i=a-b+1}^{n-b+1}\Gamma(z_iz_j^{-1},t^2z_iz_j)} \\
&\times\prod_{i=1}^{n-a} \frac{
\prod_{j=b+1}^{n+1}\Gamma(tw_i^{\pm}x_j)
\prod_{j=a-b+1}^{n-b+1} \Gamma(t^{-1}w_i^{\pm}z_j^{-1})}
{\Gamma(w_i^{\pm 2})
\prod_{j=1}^{a-b}\Gamma(tw_i^{\pm}z_j)} \\
&\times
\prod_{1\leq i<j\leq a-b} \frac{1}{\Gamma(z_iz_j^{-1},z_i^{-1}z_j,t^2z_iz_j,t^2z_iz_j)}
\\
&\times
\prod_{a-b+1\leq i<j\leq n-b+1} \frac{1}
{\Gamma(z_iz_j^{-1},z_i^{-1}z_j,t^2z_iz_j,t^{-2}z_i^{-1}z_j^{-1})},
\end{align*}
where, given $z^{(n-b)}$,
\begin{equation}\label{znb}
z_{n-b+1}:=X_b^{-1}Z^{-1}_{n-b},
\end{equation}
i.e., $z_1\cdots z_{n-b+1}=X_b^{-1}=x_1^{-1}\cdots x_b^{-1}$.
For $b=n$ this implies $z_1=X_n^{-1}=x_{n+1}$ so that
$\Delta(z^{(0)},w^{(0)},x;t)$ is given by \eqref{wantedterm}.

In the C$_n$ elliptic elliptic beta integral \eqref{Cnbeta} we
deform $\T^n$ to $C^n=C\times\cdots\times C$
with $C=C^{-1}\subset\Complex$ the usual positively oriented Jordan curve,
such that the points $t_i p^{\mu}q^{\nu}$ for $i\in\set{2n+4}$ are in the interior of $C$.
With $\T^n$ replaced by such $C^n$ the integral \eqref{Cnbeta} holds
for all $t_1,\dots,t_{2n+4}$ subject only to $t_1\cdots t_{2n+4}=pq$.
Then choosing $t_1\cdots t_{2n+2}=1$ and $t_{2n+3}t_{2n+4}=pq$
and using \eqref{Gammasymm} and $\Gamma(pq)=0$, we get
 \begin{equation*}
\int_{C^n} \prod_{j=1}^n\frac{\prod_{i=1}^{2n+2}\Gamma(t_iz_j^{\pm})}
{\Gamma(z_j^{\pm 2})}
\prod_{1\leq i<j\leq n} \frac{1}{\Gamma(z_i^{\pm}z_j^{\pm})} \d{z}=0,
\end{equation*}
where $C^n=C\times\cdots\times C$ such that $C$ has the points
$t_ip^{\mu}q^{\mu}$ for $i\in\set{2n+2}$ in its interior.
Replacing $z$ by $w$, $n$ by $n-b$, and
$t_i\to tx_{i+b}$, $t_{i+n-b+1}\to t^{-1}z_i^{-1}$ for $i\in\set{n-b+1}$,
with $z_{n-b+1}$ defined by \eqref{znb} to ensure that
\begin{multline*}
1=t_1\cdots t_{2n+2}\to x_{b+1}\cdots x_{n+1}z_1^{-1}\cdots z_{n-b+1}^{-1} \\
=x_{b+1}\cdots x_{n+1}X_b=x_1\cdots x_{n+1}=1,
\end{multline*}
yields,
\begin{equation*}
\int_{C^{n-b}}\Delta(z^{(n-b)},w^{(n-b)},x;t)\d{w^{(n-b)}}=0
\end{equation*}
for $b\in\Zet{n}$.
Here $C^{n-b}=C\times\cdots\times C$ such that $C$ has the
points $tx_{i+b}p^{\mu}q^{\nu}$ and $t^{-1}z_i^{-1}p^{\mu}q^{\nu}$ for $i\in\set{n-b+1}$
in its interior.
Obviously, we then also have
\begin{equation}\label{nul2}
\int_{\T^{n-b}}\int_{C^{n-b}}
\Delta(z^{(n-b)},w^{(n-b)},x;t)
f(z^{(n-b)},x^{(b)})
\d{w^{(n-b)}}\d{z^{(n-b)}}=0,
\end{equation}
to be compared with \eqref{nul}.

Deforming $C^{n-b}$ to $\T^{n-b}$ by successively deforming the
$1$-dimensional contours $C$ to $\T$, the poles of the integrand at
$t^{-1}z_i^{-1}$ ($t z_i$) for $i\in\set{n-b+1}$
move from the interior (exterior) of $C$ to the exterior (interior) of $\T$,
but no other poles cross the contours of integration.
The function $\Delta(z^{(n-b)},w^{(n-a)},x;t)$ has been obtained
from $\Delta(z,w,x;t)$ by computing residues corresponding to poles in the
$z$-variables, see \eqref{res1} and \eqref{recdef}.
Presently we are at a stage of the calculation that requires the
computation of residues with respect to the above-listed poles in the $w$-variables.
Amazingly, this does not lead to a further generalization of
$\Delta(z,w,x;t)$. Indeed, the following truly remarkable relation holds
for $a\in\{b,\dots,n\}$:
\begin{multline}\label{resam}
\Res_{w_{n-a}=t^{-1}z_{a-b+1}^{-1}} \frac{\Delta(z^{(n-b)},w^{(n-a)},x;t)}{w_{n-a}} \\
=-\Res_{w_{n-a}=tz_{a-b+1}} \frac{\Delta(z^{(n-b)},w^{(n-a)},x;t)}{w_{n-a}}
=\Delta(z^{(n-b)},w^{(n-a-1)},x;t).
\end{multline}
Moreover, since $\Delta(z^{(n-b)},w^{(n-a)},x;t)$ exhibits permutation symmetry in
the variables $z_{a-b+1},\dots,z_{n-b+1}$ we also have
\begin{multline}\label{ressym}
\Res_{w_{n-a}=t^{-1}z_j^{-1}} \frac{\Delta(z^{(n-b)},w^{(n-a)},x;t)}{w_{n-a}} \\
=-\Res_{w_{n-a}=tz_j} \frac{\Delta(z^{(n-b)},w^{(n-a)},x;t)}{w_{n-a}} \\
=\Delta(z^{(n-b)},w^{(n-a-1)},x;t)|_{z_{a-b+1}\leftrightarrow z_j}
\end{multline}
for $j\in\{a-b+1,\dots,n-b+1\}$.
These results imply our next lemma.
\begin{lemma}\label{lemcomp}
For $b\in\Zet{n}$ there holds
\begin{multline*}
\sum_{a=b}^n \frac{(4\pi\im)^a}{(n-a)!}
\binom{n-b}{n-a}
\biggl[\int_{\T^{n-b}}+\frac{1}{n-a+1}\sum_{c=1}^{a-b}
\int_{\T^{c-1}\times(X_b^{-1}\T)\times\T^{n-b-c}} \biggr] \\
\times \int_{\T^{n-a}}
\Delta(z^{(n-b)},w^{(n-a)},x;t)
f(z^{(n-b)},x^{(b)})
\d{w^{(n-a)}}\d{z^{(n-b)}}=0.
\end{multline*}
\end{lemma}
Before proving the lemma we complete the proof of Proposition~\ref{propnul}.

If we can show that for all $c\in\set{a-b}$ no poles of the integrand
in the $z_c$-plane
cross the contour of integration when $X_b^{-1}\T$ is inflated to $\T$
then the identity of Lemma~\ref{lemcomp} simplifies to
\begin{multline*}
\sum_{a=b}^n\frac{(4\pi\im)^a}{(n-a+1)!}
\binom{n-b}{n-a}
\int_{\T^{n-b}}\int_{\T^{n-a}}
\Delta(z^{(n-b)},w^{(n-a)},x;t) \\ \times
f(z^{(n-b)},x^{(b)})
\d{w^{(n-a)}}\d{z^{(n-b)}}=0,
\end{multline*}
where we divided out an overall factor $(n-b+1)$.
Since this is \eqref{nul} with $(k_1,\dots,k_b)=(1,\dots,b)$
we are done with the proof of Proposition~\ref{propnul} if we can show
that $\Delta(z^{(n-b)},w^{(n-a)},x;t)$ is free of poles in the annulus
$\abs{X_b}^{-1}\leq \abs{z_j}\leq 1$ for $j\in\set{a-b}$.
Since $\Delta(z^{(n-b)},w^{(n-a)},x;t)$ has permutation symmetry in the
variables $z_1,\dots,z_{a-b}$,
it is enough to check this condition for $j=a-b$.
The rest is a matter of
case-checking, where it should be noted that, since $b\leq n-1$,
$\abs{t}^{n-1}\leq \abs{t}^b<\abs{X_b}^{-1}$.
For example, the pole at $z_{a-b}=t^{-2}x_{i+b}^{-1}p^{-\mu}q^{-\nu}$
(for $i\in\set{n-b+1}$) has absolute value
$\abs{z_{a-b}}>\abs{t}^{-1}>1$ for $i\neq n-b+1$
and absolute value $\abs{z_{a-b}}>\abs{t}^{-n-2}>1$ for $i=n-b+1$.
The pole at $z_{a-b}=x_ip^{\mu}q^{\nu}$ (for $i=b+1,\dots,n+1$)
has absolute value $\abs{z_{a-b}}>\abs{t}^{-1}>1$ if $(\mu,\nu)=(0,0)$ and
$i\neq n+1$, has absolute value $\abs{z_{a-b}}=\abs{x_{n+1}}<\abs{X_b}^{-1}$ if
$(\mu,\nu)=(0,0)$ and $i=n+1$,
has absolute value $\abs{z_{a-b}}<\abs{t}^{-1}M<\abs{t}^{n-1}<\abs{X_b}^{-1}$
if $(\mu,\nu)\neq (0,0)$ and $i\neq n+1$,
and has absolute value $\abs{z_{a-b}}<M<\abs{t}^n<
\abs{X_b}^{-1}$ if $(\mu,\nu)\neq (0,0)$ and $i= n+1$.
\end{proof}

\begin{proof}[Proof of Lemma~\ref{lemcomp}]
We will show by induction on $d$ that for $b\in\Zet{n}$ and
$d\in\{b,\dots,n\}$ there holds
\begin{multline*}
L_{b,d}:=\sum_{a=b}^d \frac{(4\pi\im)^a}{(d-a)!}
\binom{n-b}{n-a}
\biggl[\int_{\T^{n-b}}+\frac{1}{n-a+1}\sum_{c=1}^{a-b}
\int_{\T^{c-1}\times (X_b^{-1}\T)\times\T^{n-b-c}} \biggr] \\
\times \int_{\T^{d-a}\times C^{n-d}_a}
\Delta(z^{(n-b)},w^{(n-a)},x;t)
f(z^{(n-b)},x^{(b)})
\d{w^{(n-a)}}\d{z^{(n-b)}}=0.
\end{multline*}
Here $C^{n-d}_a=C_a\times\cdots\times C_a$ with
$C_a$ a contour that has the points
\begin{subequations}\label{co}
\begin{align}
\label{co1}
tx_{i+b}p^{\mu}q^{\nu} &\text{ for } i\in\set{n-b+1},\\
\label{co2}
t^{-1}z_{i+a-b}^{-1}p^{\mu}q^{\nu}& \text{ for } i\in\set{n-a+1}, \\
t^{-1}z_ip^{\mu+1}q^{\nu+1} &\text{ for } i\in\set{a-b},
\label{co3}
\end{align}
\end{subequations}
in its interior.
Since for $d=b$ we recover \eqref{nul2} it suffices to establish
that $L_{b,d}=(d-b+1)L_{b,d+1}$ for $d\in\{b,\dots,n-1\}$.

Without loss of generality we may assume that $\abs{t}^{2n+1}<\abs{C_a}<
\abs{t}^{-2n-1}$. This is compatible with the fact that the points
listed in \eqref{co1} and \eqref{co2} must lie in the interior of $C_a$, and
guarantees that the points listed in \eqref{co3} lie in its interior.
When integrating over $w^{(n-a)}$ for fixed $z^{(n-b)}$ we may also assume
that $z_{a-b+1},\dots,z_{n-b+1}$ are all distinct.

We now deform $C^{n-d}_a$ to $C^{n-d-1}_a\times\T$.
The poles of $\Delta(z^{(n-b)},w^{(n-a)},x;t)$ in the $w_{n-a}$-plane
coincide with the points listed in \eqref{co} and their reciprocals.
The usual inspection of their moduli shows that under the assumption
that $\abs{x_i}<\abs{t}^{-1}$, the only poles crossing the contour are
those corresponding to $w_{n-a}=t^{-1}z_{k+a-b}^{-1}$ for $k\in\set{n-a+1}$.
Hence, by \eqref{resam} and \eqref{ressym},
\begin{align}\label{step1}
&I_{a,d}(z^{(n-b)},x;t):=\int_{\T^{d-a}\times C^{n-d}_a}
\Delta(z^{(n-b)},w^{(n-a)},x;t)\d{w^{(n-a)}} \\
&\quad=I_{a,d+1}(z^{(n-b)},x;t) \notag  \\
&\qquad +4\pi\im\sum_{k=a-b+1}^{n-b+1}\int_{\T^{d-a}\times C^{n-d-1}_{a;k}}
\Delta(z^{(n-b)},w^{(n-a-1)},x;t)|_{z_{a-b+1}\leftrightarrow z_k}
\d{w^{(n-a-1)}}. \notag
\end{align}
Here the integration variables $w_{d-a+1}$ and $w_{n-a}$
in the first integral on the right have been permuted in order
to simplify $\T^{d-a}\times C_a^{n-d-1}\times\T$ to $\T^{d-a+1}\times C_a^{n-d-1}$,
leading to $I_{a,d+1}$.
Furthermore, $C^{n-d-1}_{a;k}=C_{a;k}\times\cdots\times C_{a;k}$ with
$C_{a;k}$ a contour that has the points listed in
\eqref{co} in its interior be it that in \eqref{co2} and \eqref{co3} the respective
conditions $i\neq k$ and $i=k$ need to be added.
This latter condition is automatically satisfied if we
demand that $\abs{t}^{2n+1}<\abs{C_{a;k}}<
\abs{t}^{-2n-1}$ which also implies that
$C_{a;a-b+1}$ may be identified with $C_{a+1}$ so that
$C^{n-d-1}_{a;a-b+1}=C^{n-d-1}_{a+1}$.

Next consider
\begin{subequations}
\begin{align}\label{multfac}
J_{a,b,d}&:=
\int_{\T^{n-b}}I_{a,d}(z^{(n-b)},x;t)f(z^{(n-b)},x^{(b)})\d{z^{(n-b)}} \\
\intertext{and}
\label{ord}
K_{a,b,c,d}&:=\int_{\T^{c-1}\times (X_b^{-1}\T)\times\T^{n-b-c}}
I_{a,d}(z^{(n-b)},x;t)f(z^{(n-b)},x^{(b)})\d{z^{(n-b)}}.
\end{align}
\end{subequations}
We will compute these integrals using \eqref{step1}, starting with $J_{a,b,d}$.
By permuting the integration variables $z_k$ and $z_{a-b+1}$ for $k\in\{a-b+1,\dots,n-b\}$,
and by using the permutation symmetry of $f$,
the sum over $k$ in \eqref{step1} (with the term $k=n-b+1$ excluded)
simply gives rise to $(n-a)$ times $J_{a+1,b,d+1}$.
The last term in the sum is to be treated somewhat differently.
Recalling the definition of $z_{n-b+1}$ in \eqref{znb}
we make the variable change $z_{a-b+1}\leftrightarrow z_{n-b+1}=Z_{n-b}^{-1}X_b^{-1}$
in the integral corresponding to this remaining term. By the A$_n$ symmetry of $f$
this leaves the integrand unchanged.
Consequently,
\begin{equation}\label{step2}
J_{a,b,d}=J_{a,b,d+1}
+4\pi\im(n-a) J_{a+1,b,d+1}
+4\pi\im
K_{a+1,b,a-b+1,d+1}.
\end{equation}
We carry out exactly the same variables changes to compute $K_{a,b,c,d}$
for $c\in\set{a-b}$. Note that since $c\leq a-b$ these changes
do not interfere with the structure of the contours of \eqref{ord}.
The only notable difference will be that when permuting
$z_{a-b+1}\leftrightarrow z_{n-b+1}=Z_{n-b}^{-1}X_b^{-1}$
corresponding to the last term in the sum over $k$ in \eqref{step1},
the $z_{a-b+1}$ contour of integration will not change from $\T$
to $X_b^{-1}\T$ as before, but will remain just $\T$.
Indeed, when $z_{a-b+1}\to X_b^{-1}Z_{n-b}^{-1}$
we find $\abs{z_{a-b+1}}\to\abs{X_b^{-1}Z_{n-b}^{-1}}=1$
since $\abs{z_c}=\abs{X_b^{-1}}$.
Therefore
\begin{equation}\label{step3}
K_{a,b,c,d}=K_{a,b,c,d+1}+4\pi\im(n-a+1)K_{a+1,b,c,d+1}.
\end{equation}

{}The remainder of the proof is elementary.
By definition,
\begin{equation*}
L_{b,d}=\sum_{a=b}^d \frac{(4\pi\im)^a}{(d-a)!}
\binom{n-b}{n-a}
\biggl[J_{a,b,d}+\frac{1}{n-a+1}\sum_{c=1}^{a-b}K_{a,b,c,d}\biggr].
\end{equation*}
Substituting \eqref{step2} and \eqref{step3} and shifting the
summation index from $a$ to $a-1$ in all terms of the summand
that carry the subscript $a+1$ yields $L_{b,d}=(d-b+1)L_{b,d+1}$.
\end{proof}

\section{Alternative proof of the new A$_n$ elliptic beta integral}
\label{secaltproof}

The proof of Theorem~\ref{thmbetaA} presented in this section
adopts and refines a technique for proving elliptic beta integrals
via $q$-difference equations that was recently developed in \cite{Spiridonov04}.
We note that this method is different from the $q$-difference approach of
Gustafson \cite{Gustafson92}. The latter employs $q$-difference
equations depending solely on the ``external'' parameters in
beta integrals (like the $t_i$ and $s_i$ in \eqref{Anbeta} and \eqref{Cnbeta})
and not on the integration variables themselves. Although Gustafson's
method can also be applied to the integral of Theorem~\ref{thmbetaA}
by virtue of the fact that both sides of the integral identity satisfy
\begin{multline*}
\sum_{i=1}^{n+1}
\frac{\theta(q^{-1}Ss_i^{-1}s_{n+2}^{-1})}{\theta(Ss_{n+2}^{-1}s_{n+3}^{-1})}
\prod_{\substack{j=1 \\j\neq i}}^{n+1}
\frac{\theta(q^{-1}s_j^{-1}s_{n+3})}{\theta(s_is_j^{-1})} \\
\times
I(s_1,\dots,s_{i-1},qs_i,s_{i+1},\dots,s_{n+2},q^{-1}s_{n+3})
=I(s_1,\dots,s_{n+3}),
\end{multline*}
the proof would require a non-trivial vanishing hypothesis similar
to those formulated in \cite{vDS01}.

After these preliminary remarks we turn our attention to the
actual proof of Theorem~\ref{thmbetaA}.
We begin by noting that the kernel $\rho(z;s,t)$ defined in
\eqref{rhodef} satisfies a $q$-difference
equation involving the integration variables.
For $v=(v_1,\dots,v_n)\in\Complex^n$ and $a\in\Complex$ let
$\pi_{i,a}(v)=(v_1,\dots,v_{i-1},av_i,v_{i+1},\dots,v_n)$.
\begin{lemma}\label{qdifflemma}
We have
\begin{equation}\label{qdifference}
\rho(z;s,t)-\rho(z;s,\pi_{1,q}(t))
=\sum_{i=1}^n\bigl[g_i(z;s,t)-g_i(\pi_{i,q^{-1}}(z);s,t)\bigr],
\end{equation}
where $g_i(z;s,t)=\rho(z;s,t)f_i(z;s,t)$ and
\begin{multline}\label{gdef}
f_i(z;s,t)=
t_1z_{n+1}\theta(St_1^2)
\prod_{j=1}^n \frac{\theta(t_1z_j)}{\theta(z_jz_{n+1}^{-1},St_jz_{n+1}^{-1})}
\prod_{j=2}^n \frac{\theta(t_jz_i,q^{-1}St_jz_i^{-1})}{\theta(q^{-1}t_jz_{n+1})}
\\ \times
\prod_{j=1}^{n+3}\frac{\theta(s_jz_{n+1}^{-1})}{\theta(t_1s_j)}
\prod_{\substack{j=1 \\ j\neq i}}^n
\frac{\theta(q^{-1}z_j^{-1}z_{n+1},Sz_j^{-1}z_{n+1}^{-1})}
{\theta(z_iz_j^{-1},q^{-1}Sz_i^{-1}z_j^{-1})}.
\end{multline}
\end{lemma}

\begin{proof}
{}From \eqref{Gammaqd} and $\theta(x)=-x\,\theta(x^{-1})$ it readily follows that
\begin{equation*}
\frac{\rho(z;s,\pi_{1,q}(t))}{\rho(z;s,t)}=
\prod_{i=1}^{n+1}\frac{\theta(t_1z_i)}{\theta(St_1z_i^{-1})}
\prod_{j=1}^{n+3}\frac{\theta(St_1s_j^{-1})}{\theta(t_1s_j)}
\end{equation*}
and
\begin{multline*}
\frac{\rho(\pi_{i,q^{-1}}(z);s,t)}{\rho(z;s,t)}=
-\frac{z_i\theta(q^{-2}z_iz_{n+1}^{-1})}
{qz_{n+1}\theta(z_i^{-1}z_{n+1})}
\prod_{j=1}^n\frac{\theta(t_jz_{n+1},q^{-1}St_jz_{n+1}^{-1})}
{\theta(q^{-1}t_jz_i,St_jz_i^{-1})} \\
\times
\prod_{j=1}^{n+3}\frac{\theta(s_jz_i^{-1})}{\theta(q^{-1}s_jz_{n+1}^{-1})}
\prod_{\substack{j=1 \\ j\neq i}}^n
\frac{\theta(q^{-1}z_iz_j^{-1},q^{-1}z_jz_{n+1}^{-1},
Sz_i^{-1}z_j^{-1})}{\theta(z_i^{-1}z_j,z_j^{-1}z_{n+1},
q^{-1}Sz_j^{-1}z_{n+1}^{-1})}.
\end{multline*}
Dividing both sides of \eqref{qdifference} by $\rho(z;s;t)$ we thus obtain
the theta function identity
\begin{multline*}
1-\prod_{i=1}^{n+1}\frac{\theta(t_1z_i)}{\theta(St_1z_i^{-1})}
\prod_{j=1}^{n+3}\frac{\theta(St_1s_j^{-1})}{\theta(t_1s_j)}
=\sum_{i=1}^n \frac{t_1z_i\theta(St_1^2)}{\theta(St_1z_i^{-1})}
\prod_{j=1}^{n+3}\frac{\theta(s_jz_i^{-1})}{\theta(t_1s_j)}
\prod_{\substack{j=1 \\ j\neq i}}^{n+1}
\frac{\theta(t_1z_j)}{\theta(z_i^{-1}z_j)} \\
+t_1z_{n+1}\theta(St_1^2)
\sum_{i=1}^n
\prod_{j=1}^n \frac{\theta(t_1z_j)}{\theta(z_jz_{n+1}^{-1},St_jz_{n+1}^{-1})}
\prod_{j=2}^n\frac{\theta(t_jz_i,q^{-1}St_jz_i^{-1})}{\theta(q^{-1}t_jz_{n+1})} \\
\times
\prod_{j=1}^{n+3}\frac{\theta(s_jz_{n+1}^{-1})}{\theta(t_1s_j)}
\prod_{\substack{j=1 \\ j\neq i}}^n
\frac{\theta(q^{-1}z_j^{-1}z_{n+1},Sz_j^{-1}z_{n+1}^{-1})}
{\theta(z_iz_j^{-1},q^{-1}Sz_i^{-1}z_j^{-1})}.
\end{multline*}
Observing that the left-hand side and the first sum on the right
are independent of $t_2,\dots,t_n$ suggests that the above
identity is a linear combination of
\begin{equation}\label{theta1}
\sum_{i=1}^{n+1}
\frac{t_1z_i\theta(St_1^2)}{\theta(St_1z_i^{-1})}
\prod_{j=1}^{n+3}\frac{\theta(s_jz_i^{-1})}{\theta(t_1s_j)}
\prod_{\substack{j=1 \\ j\neq i}}^{n+1}
\frac{\theta(t_1z_j)}{\theta(z_i^{-1}z_j)} \\
=1-\prod_{i=1}^{n+1}\frac{\theta(t_1z_i)}{\theta(St_1z_i^{-1})}
\prod_{j=1}^{n+3}\frac{\theta(St_1s_j^{-1})}{\theta(t_1s_j)},
\end{equation}
and
\begin{equation}\label{theta2}
\sum_{i=1}^n
\prod_{j=2}^n\frac{\theta(t_jz_i,q^{-1}St_jz_i^{-1})}
{\theta(q^{-1}t_jz_{n+1},St_jz_{n+1}^{-1})}
\prod_{\substack{j=1 \\ j\neq i}}^n
\frac{\theta(q^{-1}z_j^{-1}z_{n+1},Sz_j^{-1}z_{n+1}^{-1})}
{\theta(z_iz_j^{-1},q^{-1}Sz_i^{-1}z_j^{-1})}=1.
\end{equation}

To prove \eqref{theta1} we first note that the conditions
$S=s_1\dots s_{n+3}$ and $z_1\cdots z_{n+1}=1$ may be replaced
by the single condition $Sz_1\cdots z_{n+1}=s_1\cdots s_{n+3}$.
Since this requires departing from the convention that $z_1\cdots z_{n+1}=1$ it
is perhaps better to state this generalization with $n$ replaced by $n-1$, i.e,
\begin{equation}\label{theta3}
\sum_{i=1}^n
\frac{az_i\theta(Ba^2)}{\theta(Baz_i^{-1})}
\prod_{j=1}^{n+2}\frac{\theta(b_jz_i^{-1})}{\theta(ab_j)}
\prod_{\substack{j=1 \\ j\neq i}}^n
\frac{\theta(az_j)}{\theta(z_i^{-1}z_j)} \\
=1-\prod_{i=1}^n\frac{\theta(az_i)}{\theta(Baz_i^{-1})}
\prod_{j=1}^{n+2}\frac{\theta(Bab_j^{-1})}{\theta(ab_j)},
\end{equation}
for $Bz_1\cdots z_n=b_1\cdots b_{n+2}$.

To prove \eqref{theta3} we bring all terms to one side
and write the resulting identity as $f(a)=0$.
Using $\theta(px)=-x^{-1}\theta(x)$ it is easily checked that
$f(pa)=f(a)$.

The function $f$ has poles at $a=b_j^{-1}p^m$ and $a=B^{-1}z_ip^m$
for $m\in\Z$, $j\in\set{n+2}$ and $i\in\set{n}$.
If we can show that the residues at these poles vanish
then $f(a)$ must be constant by Liouville's theorem.
That this constant must then be zero easily follows by taking $a=z_1^{-1}$.

By the permutation symmetry of $f$ in $z_1,\dots,z_n$ and
$b_1,\dots,b_{n+2}$, and by the periodicity of $f$ it suffices to consider
the poles at $a=b_1^{-1}$ and $a=B^{-1}z_1$.
By $\theta(x)=-x\theta(x^{-1})$ the residue of $f$ at
the latter pole is easily seen to vanish.
Equating the residue of $f$ at $a=b_1^{-1}$ to zero, and replacing
$(b_2,\dots,b_{n+2})\to(a_1,\dots,a_{n+1})$ and $Bb_1^{-1}\to A$, yields the
A$_n$ elliptic partial fraction expansion \cite[Equation (4.3)]{Rosengren02}
\begin{equation*}
\sum_{i=1}^n
\prod_{j=1}^{n+1}\frac{\theta(a_jz_i^{-1})}{\theta(Aa_j^{-1})}
\prod_{\substack{j=1 \\ j\neq i}}^n
\frac{\theta(Az_j^{-1})}{\theta(z_i^{-1}z_j)}=1
\end{equation*}
for $Az_1\cdots z_n=a_1\cdots a_{n+1}$.

The task of proving \eqref{theta2} is even simpler.
{}From \cite[Lemma 4.14]{Gustafson87} we have the elliptic
partial fraction expansion
\begin{equation*}
\sum_{i=1}^n
\prod_{j=1}^{n-1}\frac{\theta(b_jz_i,b_jz_i^{-1})}
{\theta(ab_j,a^{-1}b_j)}
\prod_{\substack{j=1 \\ j\neq i}}^n
\frac{\theta(az_j,a^{-1}z_j)}
{\theta(z_i^{-1}z_j,z_iz_j)}=1.
\end{equation*}
Making the substitutions
\begin{equation*}
z_i\to q^{-1/2}S^{1/2}z_i^{-1},\quad
b_j\to q^{-1/2}S^{1/2}t_j, \quad
a\to q^{-1/2}S^{-1/2}z_{n+1}
\end{equation*}
results in \eqref{theta2}.
\end{proof}

By integrating the identity of Lemma~\ref{qdifflemma} over
$\T^n$ and by rescaling some of the integration variables,
we obtain
\begin{equation}\label{int-eqn}
I(s,t)-I(s,\pi_{1,q}(t))=\kappa^{\A}\sum_{i=1}^n
\biggl(\int_{\T^n}-\int_{\T^{i-1}\times(q^{-1}\T)\times \T^{n-i}}\biggr)
g_i(z;s,t)\d{z},
\end{equation}
where $I(s,t)$ denotes the integral on the left of \eqref{rhoint}.
A careful inspection learns that $g_i(z;s,t)$ has poles in the $z_i$-plane at
\begin{equation*}
z_i=
\begin{cases}
Sz_j^{-1}p^{\mu}q^{\nu-1}  \\
St_jp^{-\mu-1}q^{-\nu-2+\delta_{j,1}} \\
t_j^{-1}p^{-\mu}q^{-\nu-1} \\
s_jp^{\mu}q^{\nu} \\
\end{cases} \hspace{-5mm},\quad
z_{n+1}=
\begin{cases}
Sz_j^{-1}p^{\mu}q^{\nu+1} & j\in\set{n}/\{i\} \\
St_jp^{-\mu-1}q^{-\nu}  & j\in\set{n} \\
t_j^{-1}p^{-\mu}q^{1-\nu-\delta_{j,1}}\qquad & j\in\set{n} \\
s_jp^{\mu}q^{\nu+1}  & j\in\set{n+3},
\end{cases}
\end{equation*}
where $z_{n+1}=a$ stands for $z_i^{-1}=az_1\cdots z_{i-1}z_{i+1}\cdots z_n$.
Imposing the conditions (which are stronger than \eqref{stcond})
\begin{multline}\label{stcond2}
\max\{\abs{t_1},\abs{q^{-1}t_2},\dots,\abs{q^{-1}t_n},
\abs{s_1},\dots,\abs{s_{n+3}}, \\
\abs{pS^{-1}t_1^{-1}},\dots,\abs{pS^{-1}t_n^{-1}},\abs{q^{-1}S}\}<1
\end{multline}
it follows that none of the listed poles of $g_i$ lies on the annulus
$1\leq\abs{z_i}\leq\abs{q}^{-1}$.
Consequently, when \eqref{stcond2} holds the right-hand side of
\eqref{int-eqn} vanishes and
\begin{equation}\label{IIq}
I(s,t)=I(s,\pi_{1,q}(t)).
\end{equation}

Expanding $I(s,t)$ in a Taylor series in $p$ we have
\begin{equation}\label{Taylor}
I(s,t;q,p)=\sum_{j=0}^{\infty} I_j(s,t;q)p^j,
\end{equation}
with $I_j(s,t;q)$ holomorphic in $s$ and $t$ for
\begin{equation*}
\max\{\abs{t_1},\dots,\abs{t_n},\abs{s_1},\dots,\abs{s_{n+3}}\}<1.
\end{equation*}
(The remaining conditions of \eqref{stcond} involving the parameter $p$
ensure convergence of the series \eqref{Taylor}, but, obviously,
bear no relation to the analyticity of $I_j(s,t;q)$.)

Thanks to \eqref{IIq} we have the termwise $q$-difference
$I_j(s,t;q)=I_j(s,\pi_{1,q}(t))$ when
\begin{equation*}
\max\{\abs{t_1},\abs{q^{-1}t_2},\dots,\abs{q^{-1}t_n},
\abs{s_1},\dots,\abs{s_{n+3}},\abs{q^{-1}S}\}<1.
\end{equation*}
Since this may be iterated and since the limiting point $t_1=0$ of the
sequence $t_1,t_1q,t_1q^2,\dots$ lies inside the domain of analyticity of
$I_j(s,t)$, we conclude that $I_j(s,t)$ is independent of $t_1$.
Lifting this to $I(s,t)$ and exploiting the symmetry in the $t_i$ it
follows that $I(s,t)$ is independent of $t$ for
\begin{equation*}
\max\{\abs{q^{-1}t_1},\dots,\abs{q^{-1}t_n},
\abs{s_1},\dots,\abs{s_{n+3}},
\abs{pS^{-1}t_1^{-1}},\dots,\abs{pS^{-1}t_n^{-1}},\abs{q^{-1}S}\}<1
\end{equation*}
and thus, by analytic continuation, for \eqref{stcond}.

In order to compute $I(s,t)=I(s)$ we repeat the reasoning of Section~\ref{seccons}
and note that when \eqref{stcond} is replaced by \eqref{stcondp}
then $I(s)$ is given by the integral on the left of \eqref{rhointC}.
We also know that this integral does not depend on $t$ and hence,
by \eqref{limint},
\begin{equation*}
I(s)=\lim_{\substack{t_i\to s_i^{-1}\\ \forall\,i\in\set{n}}} I(s)
=\kappa^\A \lim_{\substack{t_i\to s_i^{-1}\\ \forall\,i\in\set{n}}}
\sum_{\la} \rho_{\la}(s,t).
\end{equation*}
The expression on the right equals to $1$ thanks to the
trivial fact that for $N_i=0$ we have $\lambda=0$ and
the conditions $t_i=s_i^{-1}$ yield $\kappa^\A \rho_0(s,t)=1$.
As a final remark, we note that application of the
full residue calculus with $N_i\neq 0$ to our A$_n$ integral
results in Theorems~\ref{thmRos} and \ref{thmeB}.
Therefore, the considerations of the present section provide
an alternative proof of the corresponding sums as well.

\bibliographystyle{amsplain}

\end{document}